\font\title =cmbx10 at 12pt
\magnification =1200
\parindent=10pt
\baselineskip=18pt
\def\Bbb N{\hbox{N}}
\def\d{\hbox{d}}
\def\ui{\hbox{i}}
\def\Ui{{\rm{i}}}
\centerline{{\title 
Second-order linear differential equations with two irregular singular}}
\centerline{{\title points of rank three: the characteristic exponent}}
\bigskip
\centerline{{Wolfgang B\"uhring}}
\bigskip
\centerline{{\it Physikalisches Institut, Universit\"at Heidelberg, Philosophenweg 12,}}
 \centerline{{\it D-69120 Heidelberg, Germany}}
\bigskip

\noindent
{\bf Abstract}

For a second-order linear differential equation with two irregular singular points of rank three, multiple Laplace-type contour integral solutions are considered. An explicit formula in terms of the Stokes multipliers is derived for the characteristic exponent of the multiplicative solutions. The Stokes multipliers are represented by converging series with terms for which limit formulas as well as more detailed asymptotic expansions are available. Here certain new, recursively known coefficients enter , which are closely related to but different from the coefficients of the formal solutions at one of the irregular singular points of the differential equation. The coefficients of the formal solutions then appear as finite sums over subsets of the new coefficients. As a by-product, the leading exponential terms of the asymptotic behaviour of the late coefficients of the formal solutions are given, and this is a concrete example of the structural results obtained by Immink in a more general setting. The formulas displayed in this paper are not of merely theoretical interest, but they also are complete in the sense that they could be (and have been) implemented for computing accurate numerical values of the characteristic exponent, although the computational load is not small and increases with the rank of the singular point under consideration.
\bigskip

\bigskip

\bigskip

\noindent
{\bf 1.  Introduction}

Let us consider the differential equation
$$z^2f''+zf'-[\sum\limits_{m=1}^6 {D_mz^{-m}}+L^2+\sum\limits_{m=1}^6 {B_mz^m}]f(z)=0\eqno (1.1)$$
with the thirteen parameters $D_1,\ldots , D_6$,  $L$,  $B_1,\ldots , B_6$, which for simplicity of presentation are assumed to be real. This differential equation has two irregular singular points, each of rank 3 (if $D_6\neq 0$ and $B_6\neq 0$, respectively), at the origin and at infinity. Without loss of generality, one of the parameters except $L$ could be set equal to $1$ . We assume that $B_6$ is positive, again for simplicity of presentation, and that $L$ is not negative.

At infinity, there are formal power series solutions 
$$f_{\infty 1}^{\rm{asy}}(z)=\exp (P(z))z^{-\tau (1)}\sum\limits_{n=0}^\infty  {a_n(1)z^{-n}},\eqno(1.2)$$
$$f_{\infty 2}^{\rm{asy}}(z)=\exp (-P(z))z^{-\tau (-1)}\sum\limits_{n=0}^\infty  {a_n(-1)z^{-n}},\eqno(1.3)$$
where the various quantities are determined by
$$P(z)=p_3z^3+p_2z^2+p_1z,\eqno(1.4)$$
$$p_3={\textstyle{1 \over 3}}\sqrt {B_6},\quad\quad p_2={\textstyle{1 \over {12}}}B_5/p_3,\quad\quad p_1={\textstyle{1 \over 6}}(B_4-4p_2^2)/p_3,\eqno(1.5)$$
$$\tau (\kappa )={\textstyle{3 \over 2}}-{\textstyle{1 \over 6}}{{B_3} \over {\kappa p_3}}+{\textstyle{2 \over 3}}{{p_1p_2} \over {\kappa p_3}},\eqno(1.6)$$
$$\kappa \in \{1,-1\}\eqno(1.7)$$
and where the coefficients
$$a_n(\kappa )=a_n\eqno(1.8)$$
are known recursively by $$a_0=1\eqno(1.9)$$ and
$$\eqalign{&6\kappa p_3na_n=[-4\kappa p_2(\tau (\kappa )+n-2)+p_1^2-B_2]a_{n-1}\cr
  &+[-2\kappa p_1(\tau (\kappa )+n-{\textstyle{5 \over 2}})-B_1]a_{n-2}+(\tau (\kappa )+n-3-L)(\tau (\kappa )+n-3+L)a_{n-3}\cr
  &-\sum\limits_{m=1}^6 {D_ma_{n-m-3}},\cr}\eqno(1.10)$$
$$n=1,2,\ldots ,\quad\quad(a_{-8}=a_{-7}=\ldots =a_{-1}=0).$$
The formal solutions are asymptotic expansions as $z\to\infty$ in appropriate sectors of the complex plane.

At the origin, there are analogical formal solutions, which may be obtained from those at infinity by the simultaneous replacements $z\leftrightarrow 1/z ,\quad B_m\leftrightarrow D_m, \quad m=1,\ldots , 6$. This symmetry is our main reason for choosing just (1.1) as the standard form for a differential equation with two irregular singular points.

In the ring-shaped region $0<|z|<\infty$ we have (convergent) Floquet solutions $f^{(\omega )} (z)$ and $f^{(-\omega )} (z)$, where
$$f^{(\omega )}(z)=z^{\omega }\sum\limits_{n=-\infty }^{+\infty } {c_n^{(\omega )}z^n},\eqno(1.11)$$
which are linearly independent if $2\omega$ is not equal to an integer. Here the coefficients $c_n^{(\omega )}$ obey the recurrence relation
$$-\sum\limits_{m=1}^6 {B_mc_{n-m}^{(\omega )}}+(-L+\omega +n)(L+\omega +n)c_n^{(\omega )}-\sum\limits_{m=1}^6 {D_mc_{n+m}^{(\omega )}}=0,\eqno(1.12)$$
and we want to normalize them by choosing
$$c_0^{(\omega )}=1.\eqno(1.13)$$
The requirement that the power series converge determines, modulo 1, the possible values of the characteristic exponent (or circuit exponent or Floquet exponent) $\omega$. The problem to compute the characteristic exponent appears also and is best known  in the context of Hill's differential equation [22]. There are methods to compute the characteristic exponent numerically, which require, for $\omega=0$, the evaluation of the infinite determinant [20] associated with (1.12), or numerical integration of the differential equation along a suitable contour [21], or numerical solution of an eigen-value problem [14].

This paper developes an entirely different method for evaluating the characteristic exponent. We obtain an explicit formula in terms of quantities which are essentially the Stokes multipliers, and these are given explicitly as convergent series, the terms of which are represented by asymptotic expansions. Here certain recursively known coefficients enter which are closely related to but different from the coefficients of the formal solutions (1.2)-(1.3).

Although some of the other authors concerned with irregular singular points of rank larger than one [1, 2, 4-6, 8-10, 12, 13, 15, 18, 19] consider the general case of arbitrary rank , we here prefer to restrict our attention to rank three. This is already general enough to give an impression of what can be expected in the case of even higher rank. On the other hand, it is still simple enough so that we can, for the relevant quantities, obtain explicit expressions which are not only theoretically interesting but can also be implemented (and have been implemented) for numerical evaluation. This work may be viewed as an attempt to extend [3], which was useful for rank one, to the much more complicated case of higher rank.
\bigskip

\noindent
{\bf 
2.  Laplace contour integral solutions}

We try to apply the classical method of multiple Laplace contour integral solutions [7] and write for a solution of (1.1)
$$f(z)=z^\lambda (2\pi\ui)^{-3}\int_{C_{t_2}} {\int_{C_{t_1}} {\int_{C_s} {\exp (z^3s+z^2t_1+zt_2)v(s,t_1,t_2)}\d s\d t_1\d t_2}},\eqno(2.1)$$
where a  power factor with a still arbitrary parameter $\lambda$ has been included in view of later benefits. To derive the appropriate weight function $v(s,t_1,t_2)$ is a somewhat lengthy but not principally difficult procedure. We therefore give two lemmata stating the results of this procedure and postpone the proofs to a later section.

{\bf Lemma 1.}
{\it  The weight function  $v(s,t_1,t_2)$ has to be a solution of the partial differential equation 
$$\eqalign{&9(s^2-s_0^2)(\partial ^4v/\partial s^4)+12(st_1-s_0t_{10})(\partial ^4v/\partial s^3\partial t_1)+6(st_2-s_0t_{20})(\partial ^4v/\partial s^3\partial t_2)\cr
  &+4(t_1^2-t_{10}^2)(\partial ^4v/\partial s^2\partial t_1^2)+(t_2^2-t_{20}^2)(\partial ^4v/\partial s^2\partial t_2^2)+4(t_1t_2-t_{10}t_{20})(\partial ^4v/\partial s^2\partial t_1\partial t_2)\cr
  &+([81-6\lambda ]s+[B_3-4t_{10}t_{20}])(\partial ^3v/\partial s^3)+([52-4\lambda ]t_1+[B_2-t_{20}^2])(\partial ^3v/\partial s^2\partial t_1)\cr
  &+([25-2\lambda ]t_2+B_1)(\partial ^3v/\partial s^2\partial t_2)+([\lambda -12]^2-L^2)(\partial ^2v/\partial s^2)-D_1(\partial ^2v/\partial s\partial t_1)\cr
  &-D_2(\partial ^2v/\partial s\partial t_2)+D_3(\partial v/\partial s)+D_4(\partial v/\partial t_1)+D_5(\partial v/\partial t_2)-D_6v=0\cr}\eqno(2.2)$$
with two finite singular points at  $(s,t_1,t_2)=(\kappa s_0, \kappa t_{10},\kappa t_{20})$ where 
$$t_{20}=p_1,\quad\quad t_{10}=p_2,\quad\quad s_0=p_3,\eqno(2.3)$$
and the contours for each variable have to satisfy the condition that a certain lengthy expression, bilinear in $K$ and $v$ or their partial derivatives, have the same value at the start and the end of the contour.}

{\bf Lemma 2.}
{\it  For $\kappa\in\{1,-1\}$, there are appropriate weight functions $v=V(\kappa;s,t_1,t_2)$ which at the singular point $(s,t_1,t_2)=(\kappa s_0, \kappa t_{10},\kappa t_{20})$ have the power series expansion
$$\eqalign{&V(\kappa ;s,t_1,t_2)\cr
  &=\sum\limits_{m=0}^\infty  {\sum\limits_{n_1=0}^\infty  {\sum\limits_{n_2=0}^\infty  {A(\kappa ;m,n_1,n_2)(s-\kappa s_0)^{\mu(\kappa) +m}(t_1-\kappa t_{10})^{-\nu _1-n_1}(t_2-\kappa t_{20})^{-\nu _2-n_2}}},}\cr}\eqno(2.4)$$
where the $A$-coefficients satisfy a certain recurrence relation and the exponents are related by
$$ 3\mu (\kappa )-2\nu _1-\nu _2=\lambda+\tau (\kappa )-6.\eqno(2.5)$$
An appropriate set of coefficients is 
$$A(\kappa ;m,n_1,n_2)=\Gamma (-\mu (\kappa )-m)\Gamma (\nu _1+n_1)\Gamma (\nu _2+n_2)b(\kappa ;m,n_1,n_2),\eqno(2.6)$$
where the new coefficients $b(\kappa ;m,n_1,n_2)$ are given by the recurrence relation
$$\eqalign{&6\kappa s_0(3m-2n_1-n_2)b(\kappa ;m,n_1,n_2)\cr
  &=[(3m-2n_1-n_2+\tau(\kappa) -3)^2-L^2]b(\kappa ;m-1,n_1,n_2)\cr
  &+[-4\kappa t_{10}(3m-2n_1-n_2+\tau(\kappa) -2)+t_{20}^2-B_2]b(\kappa ;m-1,n_1-1,n_2)\cr
  &+[-2\kappa t_{20}(3m-2n_1-n_2+\tau(\kappa) -{\textstyle{5 \over 2}})-B_1]b(\kappa ;m-1,n_1,n_2-1)\cr
  &-D_1b(\kappa ;m-2,n_1-1,n_2)-D_2b(\kappa ;m-2,n_1,n_2-1)-D_3b(\kappa ;m-2,n_1,n_2)\cr
  &-D_4b(\kappa ;m-3,n_1-1,n_2)-D_5b(\kappa ;m-3,n_1,n_2-1)-D_6b(\kappa ;m-3,n_1,n_2)\cr}\eqno(2.7)$$
with the initial conditions 
$$b(\kappa ;0,0,0)=1,$$
$$b(\kappa ;m,n_1,n_2)=0\quad\hbox{if}\quad 3m-2n_1-n_2=0\quad\hbox{for}\quad (m,n_1,n_2)\ne (0,0,0),\eqno(2.8)$$
$$(b(\kappa;m,n_1,n_2)=0\quad\hbox{if}\quad m<0\quad\hbox{or}\quad n_1<0\quad\hbox{or}\quad n_2<0)$$
and with
$$b(\kappa ;0,n_1,n_2)=0\quad\hbox{for}\quad (n_1,n_2)\neq(0,0),\eqno(2.9)$$
$$b(\kappa ;m,n_1,n_2)=0\quad\hbox{for} \quad n_1+n_2>m\eqno(2.10)$$
as a consequence. In addition, there are weight functions $U(\kappa ;s,t_1,t_2)$ which are analytic in $s$ at the respective singular point, corresponding to $\mu (\kappa )=0, 1, 2$ {\rm (}without any further relation such as {\rm (2.5))}.}

In order to avoid unnecessary complications, we want to assume that the non-trivial exponent $\mu$ according to (2.5) is not equal to an integer. This can always be guarantied by a suitable choice of the still disposable parameter $\lambda$.

Since the exponents $\mu (\kappa )$, $ \nu_1$, $\nu_2$ are restricted only by (2.5) but otherwise arbitrary, there are other solutions of the partial differential equation (2.2 ) relevant as weight functions in our contour integrals (2.1). We may assume that $ \nu_1$ and $\nu_2$ are positive integers, preferentially 
$$\nu_1=\nu_2=1,\eqno(2.11)$$
but for the time being we want to keep $ \nu_1$ and $\nu_2$ in the formulas. If $ \nu_1$ and $\nu_2$ are increased by any positive integers $q_1$ and $q_2$, respectively, and  $\mu (\kappa )$ simultaneously is decreased by $(2/3)q_1+(1/3)q_2$, then (2.5) is still satisfied. We therefore have to consider, for $q_1,q_2=0, 1, 2, \ldots $, the set of solutions
$$\eqalign{&V(\kappa ;q_1,q_2)=V(\kappa ;q_1,q_2;s,t_1,t_2)\cr
  &=\sum\limits_{m=0}^\infty  {\sum\limits_{n_1=0}^m {\sum\limits_{n_2=0}^m {\Gamma (\nu _1+q_1+n_1)\Gamma (\nu _2+q_2+n_2)}}}\Gamma (-\mu (\kappa )-{\textstyle{2 \over 3}}q_1-{\textstyle{1 \over 3}}q_2-m)\cr
  &\times b(\kappa ;m,n_1,n_2)(s-\kappa s_0)^{\mu (\kappa )+{\textstyle{2 \over 3}}q_1+{\textstyle{1 \over 3}}q_2+m}(t_1-\kappa t_{10})^{-\nu _1-q_1-n_1}(t_2-\kappa t_{20})^{-\nu _2-q_2-n_2}.\cr}\eqno(2.12)$$
As indicated, we will use a short-hand notation suppressing the dependence on the variables $s,t_1,t_2$. We now have to choose appropriate contours for the integral representation (2.1). For each of the integrals over $t_1$ or $t_2$ a closed circle, traversed once in the positive sense, around the relevant singular point of the integrand is appropriate, since the pertinent exponent is an integer. For the $s$ -integral we need an infinite contour which starts somewhere at infinity where the exponential factor of the integral vanishes, surrounds one of the singular points in the positive sense, and returns (on a different sheet) to the starting point. Assuming that $s_0$ is real and positive, then, if
$$0<\arg (z)<{\textstyle{1 \over 3}}\pi\eqno(2.13)$$
the starting- and end-point at infinity has the phase $\pi /2$. If we agree that 
$$\arg (s-\kappa s_0)=0\eqno(2.14)$$
when $s$ is positive and sufficiently large, the integral of a single term of the infinite series (2.4) can be evaluated:
$$\eqalign{&z^\lambda (2\pi\ui)^{-3}\oint\limits^{(\kappa t_{20}+)} {\oint\limits^{(\kappa t_{10}+)} {\int\limits_{\kappa s_0+\Ui\infty }^{(\kappa s_0+)} {\exp (z^3s+z^2t_1+zt_2)}}}\cr
  &\times (s-\kappa s_0)^{\mu (\kappa )+m}(t_1-\kappa t_{10})^{-\nu _1-n_1}(t_2-\kappa t_{20})^{-\nu _2-n_2}\d s\d t_1\d t_2\cr
  &=\exp (\kappa s_0z^3+\kappa t_{10}z^2+\kappa t_{20}z){{z^{-\tau (\kappa )-3m+2n_1+n_2}} \over {\Gamma (-\mu (\kappa )-m)\Gamma (\nu _1+n_1)\Gamma (\nu _2+n_2)}},\cr}\eqno(2.15)$$
where, in the power of $z$ on the right-hand side, already use has been made of (2.5). As a consequence, the integral of $V(\kappa;q_1,q_2)$ yields, if the series is integrated term by term, one or the other of the formal solutions (1.2)-(1.3):
$$\eqalign{f_{\infty j}(z):=&z^\lambda (2\pi \ui)^{-3}\oint\limits^{(\kappa t_{20}+)} {\;\;\oint\limits^{(\kappa t_{10}+)} {\int\limits_{\kappa s_0+\Ui\infty }^{(\kappa s_0+)} {\exp (z^3s+z^2t_1+zt_2)}}}\cr
  &\times V(\kappa ;q_1,q_2;s,t_1,t_2)\d s\d t_1\d t_2\sim f_{\infty j}^{\rm{asy}}(z)\cr}\eqno(2.16)$$
in the sector $0<\arg (z)<{\textstyle{1 \over 3}}\pi $, where $j=1$ if $\kappa =1$ or $j=2$ if $\kappa =-1$ . Each of the solutions defined by the integral representation (2.16 ) has one of the formal solutions (1.2) or (1.3) as its asymptotic expansion as $z\rightarrow\infty$ in the indicated sector. It follows by rotation of the contour that the asymptotic expansions are theoretically valid in the larger sectors  $-{\textstyle{1 \over 6}}\pi<\arg (z)<{\textstyle{5 \over 6}}\pi $ for $j=1$ or  $-{\textstyle{1 \over 2}}\pi<\arg (z)<{\textstyle{1 \over 2}}\pi $ for $j=2$, respectively.

Looking at (2.15) we may see that all the terms for which $3m-2n_1-n_2$ is the same yield the same power of $z$. For the coefficients of the asymptotic expansions (1.2) or (1.3)  we then have the representation
$$a_n(\kappa )=\sum\limits_{(m,n_1,n_2)\in I_n} {b(\kappa ;m,n_1,n_2)},\eqno(2.17)$$
where
$$I_n=\{(m,n_1,n_2):3m-2n_1-n_2=n\}\subset \Bbb N_0\times \Bbb N_0\times \Bbb N_0,\quad n=0,1,\ldots.\eqno(2.18)$$
Because of the properties of the $b(m,n_1,n_2)$, the sum in (2.18) is finite for each finite $n$, in particular we have 
$$a_0(\kappa)=b_0(\kappa;0,0,0)=1.\eqno(2.19)$$
\bigskip

\noindent
{\bf 3.  Analytic continuation of the integrand}

Below we have to consider the integral representation (2.1) with $t_1$ and $t_2$-contours which are simple closed curves surrounding in the positive sense both the finite singular points $-t_{10}$ and  $t_{10}$ or $-t_{20}$ and  $t_{20}$, respectively, and with an $s$-contour which starts at or near $-s_0+\ui\infty$, surrounds both the finite singular points $-s_0$ and $s_0$ once in the positive sense and ends at  $s_0+\ui\infty$. We therefore need the analytic continuation of the integrand between the two singular points along this contour. With appropriate power factors $\Phi$ included for later convenience, the continuation formula reads
$$\eqalign{&\Phi (\kappa ;r_1,r_2)V(\kappa ;r_1,r_2)\cr
  &=\sum\limits_{q_1=r_1}^\infty  {\sum\limits_{q_2=r_2}^\infty  {E(-\kappa ;r_1,r_2;q_1,q_2)\Phi (-\kappa ;q_1,q_2)V(-\kappa ;q_1,q_2)}}+U(-\kappa ;r_1,r_2),\cr}\eqno(3.1)$$
where
$$\Phi (\kappa ;r_1,r_2)=(-2\kappa s_0)^{-\mu (\kappa )-{\textstyle{2 \over 3}}r_1-{\textstyle{1 \over 3}}r_2}.\eqno(3.2)$$
The effect of these power factors is that in (3.1) the total powers with non-integer exponents are powers of 
${\textstyle{1 \over 2}}[1-s/(\kappa s_0)]$ or of ${\textstyle{1 \over 2}}[1+s/(\kappa s_0)]$, respectively. We may agree that 
$\arg (1-s/s_0)=\arg(1+s/s_0)=0$ when $s$ is on the real axis between $-s_0$ and $s_0$. Also, as above, $\arg(s)=0$ when s is larger than $s_0$. Then, by analytic continuation along the contour under consideration, we have
$$\Phi (1;r_1,r_2)=(2s_0)^{-\mu (1)-{\textstyle{2 \over 3}}r_1-{\textstyle{1 \over 3}}r_2}\exp (\ui\pi (\mu (1)+{\textstyle{2 \over 3}}r_1+{\textstyle{1 \over 3}}r_2),\eqno(3.3)$$
$$\Phi (-1;r_1,r_2)=(2s_0)^{-\mu (-1)-{\textstyle{2 \over 3}}r_1-{\textstyle{1 \over 3}}r_2}.\eqno(3.4)$$
Let us rewrite (3.1) using an even more condensed notation, writing $q$ for $(q_1, q_2)$ in the parameter list of the  various functions and writing a sum over $q$  in place of a double sum over $q_1$ and $q_2$, etc. The above continuation formula then reads
$$\Phi (\kappa ;r)V(\kappa ;r)=\sum\limits_{q=r} {E(-\kappa ;r,q)\Phi (-\kappa ;q)V(-\kappa ;q)}+U(-\kappa ;r),\eqno(3.5)$$
and the second continuation formula
$$\eqalign{&U(\kappa ;r)=\Phi (-\kappa ,r)V(-\kappa ,r)\cr
  &-\sum\limits_{q=r} {\sum\limits_{p=q} {E(\kappa ;r,q)}E(-\kappa ;q,p)\Phi (-\kappa ;p)}V(-\kappa ;p)-\sum\limits_{q=r} {E(\kappa ;r,q)}U(-\kappa ;q)\cr}\eqno(3.6)$$
follows by the requirement of consistency of (3.5) for $\kappa$ replaced by $-\kappa$.
\bigskip

\noindent
{\bf
4.  Asymptotic expansions for the coefficients in the continuation formula}

We now want to determine the $E$-coefficients in the continuation formulas (3.5) , (3.6) by means of the asymptotic method of Darboux [16] applied to the variable $s$ . The left-hand side of the continuation formula (3.5) is
$$\eqalign{&({\textstyle{1 \over 2}}-{s \over {2\kappa s_0}})^{\mu (\kappa )+{\textstyle{2 \over 3}}r_1+{\textstyle{1 \over 3}}r_2}\sum\limits_{m=0}^\infty  {\sum\limits_{n_1=0}^m {\sum\limits_{n_2=0}^m {\Gamma (\nu _1+r_1+n_1)\Gamma (\nu _2+r_2+n_2)}}}\cr
  &\times \Gamma (-\mu (\kappa )-{\textstyle{2 \over 3}}r_1-{\textstyle{1 \over 3}}r_2-m)b(\kappa ;m,n_1,n_2)(-2\kappa s_0)^m\cr
  &\times ({\textstyle{1 \over 2}}-{s \over {2\kappa s_0}})^m(t_1-\kappa t_{10})^{-\nu _1-r_1-n_1}(t_2-\kappa t_{20})^{-\nu _2-r_2-n_2}.\cr}\eqno(4.1)$$
The leading singular term, when  $(s,t_1,t_2)\to(-\kappa s_0, -\kappa t_{10},-\kappa t_{20})$, on the right-hand side is
$$\eqalign{&\sum\limits_{q_1=r_1}^\infty  {\sum\limits_{q_2=r_2}^\infty  {E(-\kappa ;}}r_1,r_2;q_1,q_2)\Gamma (-\mu (-\kappa )-{\textstyle{2 \over 3}}q_1-{\textstyle{1 \over 3}}q_2)\cr
  &\times \Gamma (\nu _1+q_1)\Gamma (\nu _2+q_2)b(-\kappa ;0,0,0)\cr
  &\times ({\textstyle{1 \over 2}}+{s \over {2\kappa s_0}})^{\mu (-\kappa )+{\textstyle{2 \over 3}}q_1+{\textstyle{1 \over 3}}q_2}(t_1+\kappa t_{10})^{-\nu _1-q_1}(t_2+\kappa t_{20})^{-\nu _2-q_2}.\cr}\eqno(4.2)$$
By means of the binomial theorem in its hypergeometric-series-form
$$(1-x)^{-\alpha }=\sum\limits_{j=0}^\infty  {{{(\alpha )_j} \over {j! }}}x^j,\eqno(4.3)$$
where
$$(\alpha )_j=\alpha (\alpha +1)\ldots (\alpha +j-1)=\Gamma (\alpha +j)/\Gamma (\alpha )$$
means the Pochhammer symbol, we may expand, if $|t|$ is sufficiently large,
$$(t+\kappa t_0)^{-\nu -q}=\sum\limits_{j=0}^\infty  {{{(\nu +q)_j} \over {{j! }}}(-2\kappa t_0)^j(t-\kappa t_0)^{-\nu -q-j}}\eqno(4.4)$$
and, if $|s/(\kappa s_0)|$ is sufficiently small, 
$$({\textstyle{1 \over 2}}+{s \over {2\kappa s_0}})^{\mu (-\kappa )+{\textstyle{2 \over 3}}q_1+{\textstyle{1 \over 3}}q_2}=\sum\limits_{m=0}^\infty  {{{(-\mu (-\kappa )-{\textstyle{2 \over 3}}q_1-{\textstyle{1 \over 3}}q_2)_m} \over {m! }}}({\textstyle{1 \over 2}}-{s \over {2\kappa s_0}})^m.\eqno(4.5)$$
Then the leading singular term on the right becomes
$$\eqalign{&\sum\limits_{q_1=r_1}^\infty  {\sum\limits_{q_2=r_2}^\infty  {E(-\kappa ;}}r_1,r_2;q_1,q_2)\Gamma (-\mu (-\kappa )-{\textstyle{2 \over 3}}q_1-{\textstyle{1 \over 3}}q_2)\cr
  &\times \Gamma (\nu _1+q_1)\sum\limits_{j_1=0}^\infty  {{{(\nu _1+q_1)_{j_1}} \over {{j_1! }}}(-2\kappa t_{10})^{j_1}(t_1-\kappa t_{10})^{-\nu _1-q_1-j_1}}\cr
  &\times \Gamma (\nu _2+q_2)\sum\limits_{j_2=0}^\infty  {{{(\nu _2+q_2)_{j_2}} \over {{j_2! }}}(-2\kappa t_{20})^{j_2}(t_2-\kappa t_{20})^{-\nu _2-q_2-j_2}}\cr
  &\times b(-\kappa ;0,0,0)\sum\limits_{m=0}^\infty  {{{(-\mu (-\kappa )-{\textstyle{2 \over 3}}q_1-{\textstyle{1 \over 3}}q_2)_m} \over {m! }}}({\textstyle{1 \over 2}}-{s \over {2\kappa s_0}})^m.\cr}\eqno(4.6)$$
The coefficients of this series should agree asymptotically, as $m \to\infty$, with those of the series on the left-hand side, where the power factor in front of the series (4.1), when $s \to -\kappa s_0$, tends to $1$ and may be omitted. We therefore obtain
$$\eqalign{&\sum\limits_{n_1=0}^m {\sum\limits_{n_2=0}^m {\Gamma (\nu _1+r_1+n_1)\Gamma (\nu _2+r_2+n_2)}}\Gamma (-\mu (\kappa )-{\textstyle{2 \over 3}}r_1-{\textstyle{1 \over 3}}r_2-m)\cr
  &\times b(\kappa ;m,n_1,n_2)(-2\kappa s_0)^m(t_1-\kappa t_{10})^{-\nu _1-r_1-n_1}(t_2-\kappa t_{20})^{-\nu _2-r_2-n_2}\cr
  &\approx \sum\limits_{q_1=r_1}^\infty  {\sum\limits_{q_2=r_2}^\infty  {E(-\kappa ;}}r_1,r_2;q_1,q_2)\Gamma (-\mu (-\kappa )-{\textstyle{2 \over 3}}q_1-{\textstyle{1 \over 3}}q_2)\cr
  &\times \Gamma (\nu _1+q_1)\sum\limits_{j_1=0}^\infty  {{{(\nu _1+q_1)_{j_1}} \over {{j_1! }}}(-2\kappa t_{10})^{j_1}(t_1-\kappa t_{10})^{-\nu _1-q_1-j_1}}\cr
  &\times \Gamma (\nu _2+q_2)\sum\limits_{j_2=0}^\infty  {{{(\nu _2+q_2)_{j_2}} \over {{j_2! }}}(-2\kappa t_{20})^{j_2}(t_2-\kappa t_{20})^{-\nu _2-q_2-j_2}}\cr
  &\times b(-\kappa ;0,0,0){{(-\mu (-\kappa )-{\textstyle{2 \over 3}}q_1-{\textstyle{1 \over 3}}q_2)_m} \over {m! }},\cr}\eqno(4.7)$$
which holds asymptotically as $m\to\infty$. On both sides of this asymptotic equation a double power series in the same variables appears, so the coefficients of the corresponding terms must be equal. This yields, for each set $r_1,r_2,n_1,n_2$,
$$\eqalign{&\Gamma (-\mu (\kappa )-{\textstyle{2 \over 3}}r_1-{\textstyle{1 \over 3}}r_2-m)(-2\kappa s_0)^mb(\kappa ;m,n_1,n_2)\cr
  &\approx \sum\limits_{q_1=r_1}^{r_1+n_1} {\sum\limits_{q_2=r_2}^{r_2+n_2} {E(-\kappa ;r_1,r_2;q_1,q_2)}}{{\Gamma (-\mu (-\kappa )-{\textstyle{2 \over 3}}q_1-{\textstyle{1 \over 3}}q_2+m)} \over {m! }}\cr
  &\times \left( {{1 \over {j_1! }}(-2\kappa t_{10})^{j_1}} \right)_{j_1=r_1+n_1-q_1}\left( {{{1} \over {j_2! }}(-2\kappa t_{20})^{j_2}} \right)_{j_2=r_2+n_2-q_2} b(-\kappa ;0,0,0)\cr}\eqno(4.8)$$
or, if we introduce new indices of summation and make use of the reflection formula of the gamma function, 
$$\eqalign{&-\pi [\Gamma (1+m)]^{-1}(2\kappa s_0)^mb(\kappa ;m,n_1,n_2)\cr
  &\approx \sum\limits_{p_1=0}^{n_1} {\sum\limits_{p_2=0}^{n_2} {E(-\kappa ;r_1,r_2;r_1+p_1,r_2+p_2)}}{{\Gamma (-\mu (-\kappa )-{\textstyle{2 \over 3}}r_1-{\textstyle{1 \over 3}}r_2-{\textstyle{2 \over 3}}p_1-{\textstyle{1 \over 3}}p_2+m)} \over {\Gamma (1+m)}}\cr
  &\times {{\Gamma (1+\mu (\kappa )+{\textstyle{2 \over 3}}r_1+{\textstyle{1 \over 3}}r_2+m)} \over {\Gamma (1+m)}}\sin (\pi [\mu (\kappa )+{\textstyle{2 \over 3}}r_1+{\textstyle{1 \over 3}}r_2])\cr
  &\times \left( {{1 \over {j_1! }}(-2\kappa t_{10})^{j_1}} \right)_{j_1=n_1-p_1}\left( {{{1} \over {j_2! }}(-2\kappa t_{20})^{j_2}} \right)_{j_2=n_2-p_2} b(-\kappa ;0,0,0).\cr}\eqno(4.9)$$
Now the left-hand side is independent of $r_1$ and $r_2$ and so is the product of the two ratios of gamma functions, when $m\to\infty$, on the right. Therefore 
$E(-\kappa  ;r_1,r_2;r_1+p_1,r_2+p_2)\sin (\pi [\mu (\kappa )+{\textstyle{2 \over 3}}r_1+{\textstyle{1 \over 3}}r_2])$ must be independent of  $r_1$ and $r_2$ too. This proves

{\bf
Lemma 3.} 
$$E(-\kappa ;r_1,r_2;r_1+p_1,r_2+p_2)={{\sin (\pi \mu (\kappa ))} \over {\sin (\pi [\mu (\kappa )+{\textstyle{2 \over 3}}r_1+{\textstyle{1 \over 3}}r_2])}}E(-\kappa ;0,0;p_1,p_2).\eqno(4.10)$$

After Lemma 3 is available, we need to determine the $E$-coefficients for $r_1=r_2=0$ only, and it is advantageous to introduce the closely related coefficients
$$e(-\kappa ;p_1,p_2)=\sin (\pi \mu (\kappa ))E(-\kappa ;0,0;p_1,p_2),\eqno(4.11)$$
which are independent of $\lambda$. Eq (4.9), with $r_1=r_2=0$ and the factor $b(-\kappa;0,0,0)$, which is equal to $1$, omitted,  then becomes
$$\eqalign{&{{-\pi } \over {\Gamma (1+\mu (\kappa )+m)}}(2\kappa s_0)^mb(\kappa ;m,n_1,n_2)\cr
  &\approx \sum\limits_{p_1=0}^{n_1} {\sum\limits_{p_2=0}^{n_2} {e(-\kappa ;p_1,p_2)}}{{\Gamma (-\mu (-\kappa )-{\textstyle{2 \over 3}}p_1-{\textstyle{1 \over 3}}p_2+m)} \over {\Gamma (1+m)}}\cr
  &\times \left( {{1 \over {j_1! }}(-2\kappa t_{10})^{j_1}} \right)_{j_1=n_1-p_1}\left( {{1 \over {j_2! }}(-2\kappa t_{20})^{j_2}} \right)_{j_2=n_2-p_2}.\cr}\eqno(4.12)$$
 We then may solve (4.12) for the  $e$-coefficient with $p_1=n_1,  p_2=n_2$ and obtain the following asymptotic formula in terms of a $b$-coefficient and the earlier $e$-coefficients,
$$\eqalign{&e(-\kappa ;n_1,n_2)\approx -\pi (2\kappa s_0)^m\cr
  &\times {{\Gamma (1+m)} \over {\Gamma (1+\mu (\kappa )+m)\Gamma (-\mu (-\kappa )-{\textstyle{2 \over 3}}n_1-{\textstyle{1 \over 3}}n_2+m)}}b(\kappa ;m,n_1,n_2)\cr
  &-\mathop {\sum\limits_{p_1=0}^{n_1} {\sum\limits_{p_2=0}^{n_2} {}}}\limits_{(p_1,p_2)\ne (n_1,n_2)}e(-\kappa ;p_1,p_2){{\Gamma (-\mu (-\kappa )-{\textstyle{2 \over 3}}p_1-{\textstyle{1 \over 3}}p_2+m)} \over {\Gamma (-\mu (-\kappa )-{\textstyle{2 \over 3}}n_1-{\textstyle{1 \over 3}}n_2+m)}}\cr
  &\times \left( {{1 \over {j_1! }}(-2\kappa t_{10})^{j_1}} \right)_{j_1=n_1-p_1}\left( {{1 \over {j_2! }}(-2\kappa t_{20})^{j_2}} \right)_{j_2=n_2-p_2}.\cr}\eqno(4.13)$$
By repeated application of this formula, all the $e$-coefficients on the right-hand side may be eliminated, and this yields the remarkable explicit limit formula
$$\eqalign{&e(-\kappa ;n_1,n_2)=-\pi \cr
  &\times \mathop {\lim }\limits_{m\to \infty }{{\Gamma (1+m)} \over {\Gamma (1+\mu (\kappa )+m)\Gamma (-\mu (-\kappa )-{\textstyle{2 \over 3}}n_1-{\textstyle{1 \over 3}}n_2+m)}}\cr
  &\times (2\kappa s_0)^m\sum\limits_{j_1=0}^{n_1} {\sum\limits_{j_2=0}^{n_2} {b(\kappa ;m,j_1,j_2)}}{{(2\kappa t_{10})^{n_1-j_1}} \over {(n_1-j_1)! }}{{(2\kappa t_{20})^{n_2-j_2}} \over {(n_2-j_2)! }}.\cr}\eqno(4.14)$$
The proof of this formula will be given below in Section 7.

Although approximate asymptotic equations or limit formulas such as (4.13) or (4.14) are interesting from a theoretical point of view, we finally need more, namely a detailed asymptotic expansion suitable for accurate numerical evaluation. This can be obtained, on the basis of Sch\"afke and Schmidt [17], essentially in the same way as above, apart from the following two refinements: The power factor in front of (4.1) can no longer be omitted, and we have to include a finite number of singular terms rather than the leading one, (4.2), alone.  The result of this procedure, which will be derived in more detail below, is

{\bf Theorem 1.} {\it The $E$-coefficients, or $e$-coefficients according to {\rm (4.11),} in the continuation formula {\rm (3.1)} or {\rm (3.5)} are
$$\eqalign{&e(-\kappa ;n_1,n_2)=-\pi \cr
  &\times \left[ {{{\Gamma (1+m)} \over {\Gamma (1+\mu (\kappa )+m)\Gamma (-\mu (-\kappa )-{\textstyle{2 \over 3}}n_1-{\textstyle{1 \over 3}}n_2+m)}}} \right.(2\kappa s_0)^mb(\kappa ;m,n_1,n_2)\cr
  &-\mathop {\sum\limits_{q_1=0}^{n_1} {\sum\limits_{q_2=0}^{n_2} {}}}\limits_{(q_1,q_2)\ne (n_1,n_2)}e(-\kappa ;q_1,q_2){{\Gamma (-\mu (-\kappa )-{\textstyle{2 \over 3}}q_1-{\textstyle{1 \over 3}}q_2+m)} \over {\Gamma (-\mu (-\kappa )-{\textstyle{2 \over 3}}n_1-{\textstyle{1 \over 3}}n_2+m)}}\cr
  &\times [1+\sum\limits_{k=1}^K {\sum\limits_{l_1=0}^k {\sum\limits_{l_2=0}^k {{{(1+\mu (-\kappa )+{\textstyle{2 \over 3}}q_1+{\textstyle{1 \over 3}}q_2)_k} \over {(1+\mu (-\kappa )+{\textstyle{2 \over 3}}q_1+{\textstyle{1 \over 3}}q_2-m)_k}}H(-\kappa ;k,l_1,l_2;q_1,q_2)}}}+O(m^{-K-1})]\cr
  &\times \left. {\left( {{1 \over {j_1! }}(-2\kappa t_{10})^{j_1}} \right)_{j_1+q_1+l_1=n_1}\left( {{1 \over {j_2! }}(-2\kappa t_{20})^{j_2}} \right)_{j_2+q_2+l_2=n_2}} \right]\cr
  &\times \left[ {1+\sum\limits_{k=1}^K {{{(1+\mu (-\kappa )+{\textstyle{2 \over 3}}n_1+{\textstyle{1 \over 3}}n_2)_k} \over {(1+\mu (-\kappa )+{\textstyle{2 \over 3}}n_1+{\textstyle{1 \over 3}}n_2-m)_k}}H(-\kappa ;k,0,0;n_1,n_2)}+O(m^{-K-1})} \right]^{-1},\cr}\eqno(4.15)$$
where}
$$H(-\kappa ;k,l_1,l_2;q_1,q_2)=\sum\limits_{j=l_1+l_2}^k {{{(\mu (\kappa ))_{k-j}} \over {(k-j)! (1+\mu (-\kappa )+{\textstyle{2 \over 3}}q_1+{\textstyle{1 \over 3}}q_2)_j}}}(-2\kappa s_0)^jb(-\kappa ;j,l_1,l_2).\eqno(4.16)$$
With a suitable choice of $m$ and $K$, Theorem 1 may be used to compute accurate values of the $e$-coefficients, beginning with
$$\eqalign{&e(-\kappa ;0,0)=-\pi {{\Gamma (1+m)} \over {\Gamma (1+\mu (\kappa )+m)\Gamma (-\mu (-\kappa )+m)}}(2\kappa s_0)^mb(\kappa ;m,0,0)\cr
  &\times \left[ {1+\sum\limits_{k=1}^K {{{(1+\mu (-\kappa ))_k} \over {(1+\mu (-\kappa )-m)_k}}H(-\kappa ;k,0,0;0,0)}+O(m^{-K-1})} \right]^{-1}.\cr}\eqno(4.17)$$
While the $e$-coefficients are independent of $\lambda$, the approximate values computed for any finite $m$ do depend on it. Thus $\lambda$ here plays the role of a computational parameter which could be adjusted for optimal accuracy. The best choice from this point of view , in particular when $L$ is not small, is $\lambda=L$
or  $\lambda=-L$, according to the discussion in Section 7 below.

\bigskip

\noindent
{\bf 5.  Multiplicative solutions and the Floquet exponent }

We now want to construct a linear combination of the solutions at infinity,
$$f_p(z)=\alpha f_{\infty 1}(z)+\beta f_{\infty 2}(z),\eqno(5.1)$$
which is a multiplicative solution such that, after analytic continuation along a sufficiently large circle around the origin traversed once in the negative sense, this solution remains the same apart from multiplication by a constant factor, that is
$$f_{p}({\rm{e}}^{-2\pi\Ui}z)=p f_{p}(z).\eqno(5.2)$$
This solution is proportional to one of the Floquet solutions introduced above in the introduction.

We use the integral representation with contours which surround both the finite singular points as introduced above in Section 3. This integral is equal to the sum of the two integrals with contours, also considered above in Section 2, which surround only one of these singular points. If the contour is kept fixed, then the circle in the $z$-plane traversed once in the negative sense corresponds to a circle, around the origin and with radius greater than $s_0$, in the $s$-plane traversed three times in the positive sense, since then $z^3s$ in the exponential part of the integrand does not change.

Let us consider the integral representation with $V(1;0)$ in the integrand and with the above phase conventions. From (3.5) we have
$$\Phi (1;0)V(1;0)=\sum\limits_{q=0} {E(-1;0,q)\Phi (-1;q)V(-1;q)}+U(-1;0).\eqno(5.3)$$
The integral then yields a function 
$$f^{\rm{I}}(z):=\Phi (1;0)f_{\infty 1}(z)+\sum\limits_{q=0} {E(-1;0,q)\Phi (-1;q)}f_{\infty 2}(z),\eqno(5.4)$$
where the two terms come from the two singular points, which contribute via the left- or right-hand side of (5.3), respectively. Let us now consider, in the $s$-plane, a simple closed loop, homotopic to the circle mentioned above, consisting of small circles around the singular points $s_0$ and $-s_0$ and straight line segments along the real axis between these singular points. Let us start with the continuation formula (5.3) in a neighbourhood of the origin, where the left- as well as the right-hand side of (5.3) are valid, and see what happens when we follow the loop in the positive sense. We shall always refer to the above phase conventions and display any additional phases explicitly. Traversing the circular part around $s_0$ multiplies the function $V(1;0)$ by the phase factor 
$\exp (2\pi\ui\mu (1 ))$, so that we then have 
$$\eqalign{&\exp (2\pi\ui\mu (1))\Phi (1;0)V(1;0)\cr
  &=\exp (2\pi\ui\mu (1))\{\sum\limits_{q=0} {E(-1;0,q)\Phi (-1;q)V(-1;q)}+U(-1;0)\},\cr}\eqno(5.5)$$
where the right-hand side is the analytic continuation of the left by means of (5.3). Next we have to traverse the circular part around  $-s_0$, which multiplies the function $V(-1;q)$ by the phase factor $\exp (2\pi\ui[\mu (-1 )+\varphi (q)])$, where
$$\varphi (q)={\textstyle{2 \over 3}}q_1+{\textstyle{1 \over 3}}q_2,\eqno(5.6)$$
 so that we have on the  straight line segment from  $-s_0$ to  $s_0$ , after the loop has been traversed once,
$$\eqalign{&\exp (2\pi\ui[\mu (1)+\mu (-1)])\sum\limits_{q=0} {E(-1;0,q)\exp (2\pi \ui\varphi (q))\Phi (-1;q)V(-1;q)}\cr
  &+\exp (2\pi\ui\mu (1))U(-1;0)\cr}$$
$$\eqalign{&=\exp (2\pi\ui\mu (1))\{\Phi (1;0)V(1;0)+\sum\limits_{q=0} {\{\exp (2\pi \ui[\mu (-1)+\varphi (q)])-1\}U(1;q)}\cr
  &+\sum\limits_{q=0} {\sum\limits_{p=q} {E(-1;0,q)\{\exp (2\pi\ui[\mu (-1)+\varphi (q)])-1\}E(1;q,p)\Phi (1;p)V(1;p)}}\},\cr}\eqno(5.7)$$
where again the right-hand side is the analytic continuation of the left. It is evident that following the loop further leads to increasingly lengthier and more complicated formulas, not suitable for being fully displayed in this paper. We therefore want to stop here for a moment and consider the integral representation with the integrand obtained after the loop in the $s$-plane has been traversed only once (rather than three times, as finally needed). A representative example of the terms in (5.7) then is
$$\sum\limits_{q=0} {\sum\limits_{p=q} {E(-1;0,q)\{\exp (2\pi\ui[\mu (-1)+\varphi (q)])-1\}E(1;q,p)\Phi (1;p)V(1;p)}}.\eqno(5.8)$$
The integral of $V(1;p)$ yields $f_{\infty1}(z)$, independent of $p$, so that the sum over $p$ can now be performed and , because of Lemma 3, the multiple sum reduces to a  product of single sums. After integration we therefore obtain for (5.8)
$$\eqalign{&\{2\ui\sin (\pi \mu (-1))\exp(\ui\pi \mu (-1))\sum\limits_{q=0} {E(-1;0,q)(2s_0)^{-\varphi (q)}\exp (2\pi\ui\varphi (q))}\cr
  &\times \sum\limits_{\tilde p=0} {E(1;0,\tilde p)(2s_0)^{-\varphi (\tilde p)}\exp (\ui\pi \varphi (\tilde p))}\}\Phi (1;0)f_{\infty 1}(z).\cr}\eqno(5.9)$$
It is now convenient to introduce the phase factor 
$$\eta :=\exp ({\textstyle{1 \over 3}}\pi\ui)={\textstyle{1 \over 2}}(1+\ui\sqrt 3)\eqno(5.10)$$
which satisfies
$$\eta ^6=1,\quad\quad 1+\eta ^2+\eta ^4=0.\eqno(5.11)$$
Also, for the sums occurring here and below we may introduce the "Stokes multipliers"
$$\eqalign{&\sigma _n(1):=\sin (\pi \mu (1))\sum\limits_{q=0} {E(-1;0,q)(2s_0)^{-\varphi (q)}\exp (2n\pi\ui\varphi (q))},\cr
  &\sigma _n(-1):=\sin (\pi \mu (-1))\sum\limits_{q=0} {E(1;0,q)(2s_0)^{-\varphi (q)}\exp ((2n+1)\pi\ui\varphi (q))},\cr}\eqno(5.12)$$
which become 
$$\eqalign{&\sigma _0(1)=S_0(1)+(2s_0)^{-1/3}S_1(1)+(2s_0)^{-2/3}S_2(1),\cr
  &\sigma _0(-1)=S_0(-1)+\eta (2s_0)^{-1/3}S_1(-1)+\eta ^2(2s_0)^{-2/3}S_2(-1),\cr
  &\sigma _1(1)=S_0(1)+\eta ^2(2s_0)^{-1/3}S_1(1)+\eta ^4(2s_0)^{-2/3}S_2(1),\cr
  &\sigma _1(-1)=S_0(-1)+\eta ^3(2s_0)^{-1/3}S_1(-1)+(2s_0)^{-2/3}S_2(-1),\cr
  &\sigma _2(1)=S_0(1)+\eta ^4(2s_0)^{-1/3}S_1(1)+\eta ^2(2s_0)^{-2/3}S_2(1),\cr
  &\sigma _2(-1)=S_0(-1)+\eta ^5(2s_0)^{-1/3}S_1(-1)+\eta ^4(2s_0)^{-2/3}S_2(-1)\cr}\eqno(5.13)$$
in terms of the real partial sums
$$\eqalign{S_0(\kappa ):=&\sum\limits_{l=0}^\infty  {(2\kappa s_0)^{-l}\sum\limits_{(n_1,n_2)\in J_{3l}} {e(-\kappa ;n_1,n_2)}}=e(-\kappa ;0,0)+\ldots ,\cr
  S_1(\kappa ):=&\sum\limits_{l=0}^\infty  {(2\kappa s_0)^{-l}\sum\limits_{(n_1,n_2)\in J_{3l+1}} {e(-\kappa ;n_1,n_2)}}=e(-\kappa ;0,1)+\ldots ,\cr
  S_2(\kappa ):=&\sum\limits_{l=0}^\infty  {(2\kappa s_0)^{-l}\sum\limits_{(n_1,n_2)\in J_{3l+2}} {e(-\kappa ;n_1,n_2)}}=e(-\kappa ;1,0)+e(-\kappa ;0,2)+\ldots ,\cr}\eqno(5.14)$$
where
$$J_l=\{(n_1,n_2):2n_1+n_2=l\}\subset\Bbb N_0\times\Bbb N_0,\quad l=0,1,\ldots .\eqno(5.15)$$
The integral of the representative term (5.9) then becomes
$${{2\ui\exp (\ui\pi \mu (-1))} \over {\sin (\pi \mu (1))}}\sigma _1(1)\sigma _0(-1)\Phi (1;0)f_{\infty 1}(z).\eqno(5.16)$$
In total, we have
$$\eqalign{&\exp (2\pi\ui\lambda /3)f^{{\rm{I}}}({\rm{e}}^{-2\pi\Ui/3}z)={{\exp (2\pi \ui[\mu (1)+\mu (-1)])} \over {\sin (\pi \mu (1))}}\sigma _1(1)\Phi (-1;0)f_{\infty 2}(z)\cr
  &+\exp (2\pi\ui\mu (1)\{1+{{2\ui\exp (\pi\ui\mu (-1))} \over {\sin (\pi \mu (1))}}\sigma _1(1)\sigma _0(-1)\}\Phi (1;0)f_{\infty 1}(z),\cr}\eqno(5.17)$$
where
$$f^{{\rm{I}}}(z)=\Phi (1;0)f_{\infty 1}(z)+{1 \over {\sin (\pi \mu (1))}}\sigma _0(1)\Phi (-1;0)f_{\infty 2}(z).\eqno(5.18)$$
So far we have traversed the $s$-loop once and obtained on it the analytic continuation of the integrand, but we have to traverse it three times. This yields
$$f^{\rm I}({\rm{e}}^{-2\pi\Ui}z)=S_{11}\Phi (1;0)f_{\infty 1}(z)+S_{12}\Phi (-1;0)f_{\infty 2}(z)\eqno(5.19)$$
with lengthy expressions for $S_{11}$ and $S_{12}$.

In a similar way, the integral representation with $V(-1;0)$ yields
$$f^{\rm{II}}(z):=\sum\limits_{q=0} {E(1;0,q)\Phi (1;q)}f_{\infty 1}(z)+\Phi (-1;0)f_{\infty 2}(z)\eqno(5.20)$$
or
$$f^{\rm{II}}(z)={1 \over {\sin (\pi \mu (-1))}}\sigma _0(-1)\Phi (1;0)f_{\infty 1}(z)+\Phi (-1;0)f_{\infty 2}(z)\eqno(5.21)$$
and leads, after the loop has been traversed once, to
$$\eqalign{&\exp (2\pi\ui\lambda /3)f^{{\rm{II}}}({\rm{e}}^{-2\pi\Ui/3}z)=\{2\ui\exp (\pi\ui\mu (-1))\sigma _0(-1)+{{\exp (2\pi\ui\mu (1))} \over {\sin (\pi \mu (-1))}}\sigma _1(-1)\cr
  &-{{4\exp (\pi\ui[\mu (1)+\mu (-1)])} \over {\sin (\pi \mu (-1))}}\sigma _1(-1)\sigma _1(1)\sigma _0(-1)\}\Phi (1;0)f_{\infty 1}(z)\cr
  &+\exp (2\pi\ui\mu (-1))\{1+{{2\ui\exp (\pi\ui\mu (1))} \over {\sin (\pi \mu (-1))}}\sigma _1(-1)\sigma _1(1)\}\Phi (-1;0)f_{\infty 2}(z)\cr}\eqno(5.22)$$
and, after three times,
$$f^{\rm{ II}}({\rm{e}}^{-2\pi\Ui}z)=S_{21}\Phi (1;0)f_{\infty 1}(z)+S_{22}\Phi (-1;0)f_{\infty 2}(z),\eqno(5.23)$$
where again the expressions for $S_{21}$ and  $S_{22}$ are too lengthy to be displayed here.

What we really want to obtain are the circuit relations for $f_{\infty 1}$ and $f_{\infty 2}$, that is
$$\eqalign{&f_{\infty 1}({\rm{e}}^{-2\pi\Ui}z)=T_{11}f_{\infty 1}(z)+T_{12}f_{\infty 2}(z),\cr
  &f_{\infty 2}({\rm{e}}^{-2\pi\Ui}z)=T_{21}f_{\infty 1}(z)+T_{22}f_{\infty 2}(z),\cr}\eqno(5.24)$$
where
$$\eqalign{&T_{11}=[S_{11}-{{\sigma _0(1)} \over {\sin (\pi \mu (1))}}S_{21}]/[1-{{\sigma _0(-1)\sigma _0(1)} \over {\sin (\pi \mu (-1))\sin (\pi \mu (1))}}],\cr
  &T_{12}=[\Phi (-1;0)/\Phi (1;0)][S_{12}-{{\sigma _0(1)} \over {\sin (\pi \mu (1))}}S_{22}]/[1-{{\sigma _0(-1)\sigma _0(1)} \over {\sin (\pi \mu (-1))\sin (\pi \mu (1))}}],\cr
  &T_{21}=[\Phi (1;0)/\Phi (-1;0)][S_{21}-{{\sigma _0(-1)} \over {\sin (\pi \mu (-1))}}S_{11}]/[1-{{\sigma _0(-1)\sigma _0(1)} \over {\sin (\pi \mu (-1))\sin (\pi \mu (1))}}],\cr
  &T_{22}=[S_{22}-{{\sigma _0(-1)} \over {\sin (\pi \mu (-1))}}S_{12}]/[1-{{\sigma _0(-1)\sigma _0(1)} \over {\sin (\pi \mu (-1))\sin (\pi \mu (1))}}].\cr}\eqno(5.25)$$
It turns out that each numerator contains a common factor which compensates the denominator. The result then is

{\bf Theorem 2.}  {\it The coefficients in the circuit relations {\rm (5.24)} are}
$$\eqalign{&T_{11}=\exp (2\pi\ui\tau (1))\cr 
&+4\exp {({\textstyle{4 \over 3}}\pi\ui\tau (-1))}\sigma _0(-1)\sigma _0(1)+4\exp {({\textstyle{4 \over 3}}\pi\ui\tau (1))}\sigma _1(1)\sigma _0(-1)+4\sigma _2(1)\sigma _0(-1)\cr
  &+4\exp {({\textstyle{4 \over 3}}\pi\ui\tau (1))}\sigma _2(1)\sigma _1(-1)+4\exp {({\textstyle{4 \over 3}}\pi\ui\tau (1))}\sigma _2(-1)\sigma _0(1)+4\sigma _1(-1)\sigma _0(1)\cr
  &+16{\exp }({\textstyle{2 \over 3}}\pi\ui\tau (1))\sigma _2(-1)\sigma _2(1)\sigma _1(-1)\sigma _0(1)+16\exp ({\textstyle{2 \over 3}}\pi\ui\tau (1))\sigma _2(-1)\sigma _1(1)\sigma _0(-1)\sigma _0(1)\cr
  &+16{\exp }({\textstyle{2 \over 3}}\pi\ui\tau (1))\sigma _2(1)\sigma _1(-1)\sigma _1(1)\sigma _0(-1)+16\exp ({\textstyle{2 \over 3}}\pi\ui\tau (-1))\sigma _2(-1)\sigma _2(1)\sigma _0(-1)\sigma _0(1)\cr
  &+16{\exp }({\textstyle{2 \over 3}}\pi\ui\tau (-1))\sigma _1(-1)\sigma _1(1)\sigma _0(-1)\sigma _0(1)\cr
  &+64\sigma _2(-1)\sigma _2(1)\sigma _1(-1)\sigma _1(1)\sigma _0(-1)\sigma _0(1),\cr}\eqno(5.26a)$$
$$\eqalign{&T_{12}=(2s_0)^{(1/3)[\tau (1)-\tau (-1)]}\{2\ui\exp ({\textstyle{2 \over 3}}\pi\ui \tau (-1))\sigma _2(1)\cr
  &+2\ui\exp ({\textstyle{2 \over 3}}\pi\ui \tau (1))\sigma _1(1)+2\ui\exp (2\pi\ui\tau (-1))\sigma _0(1)\cr
  &+8\ui\exp ({\textstyle{4 \over 3}}\pi\ui \tau (-1))\sigma _2(-1)\sigma _2(1)\sigma _0(1)+8\ui\sigma _2(-1)\sigma _1(1)\sigma _0(1)\cr
  &+8\ui\exp ({\textstyle{4 \over 3}}\pi\ui \tau (-1))\sigma _1(-1)\sigma _1(1)\sigma _0(1)+8\ui\sigma _2(1)\sigma _1(-1)\sigma _1(1)\cr
  &+32\ui\exp ({\textstyle{2 \over 3}}\pi\ui \tau (-1))\sigma _2(-1)\sigma _2(1)\sigma _1(-1)\sigma _1(1)\sigma _0(1)\},\cr}\eqno(5.26b)$$
$$\eqalign{&T_{21}=(2s_0)^{(1/3)[\tau (-1)-\tau (1)]}\{-2\ui\exp ({\textstyle{4 \over 3}}\pi\ui \tau (-1))\sigma _0(-1)\cr
  &-2\ui\exp ({\textstyle{4 \over 3}}\pi\ui \tau (1))\sigma _2(-1)-2\ui\sigma _1(-1)\cr
  &-8\ui\exp ({\textstyle{2 \over 3}}\pi\ui \tau (1))\sigma _2(-1)\sigma _2(1)\sigma _1(-1)\cr
  &-8\ui\exp ({\textstyle{2 \over 3}}\pi\ui \tau (-1))\sigma _2(-1)\sigma _2(1)\sigma _0(-1)\cr
  &-8\ui\exp ({\textstyle{2 \over 3}}\pi\ui \tau (-1))\sigma _1(-1)\sigma _1(1)\sigma _0(-1)\cr
  &-8\ui\exp ({\textstyle{2 \over 3}}\pi\ui \tau (1))\sigma _2(-1)\sigma _1(1)\sigma _0(-1)\cr
  &-32\ui\sigma _2(-1)\sigma _2(1)\sigma _1(-1)\sigma _1(1)\sigma _0(-1)\},\cr}\eqno(5.26c)$$
$$\eqalign{&T_{22}=\exp (2\pi\ui\tau (-1))\cr
  &+4{\exp }({\textstyle{4 \over 3}}\pi\ui\tau (-1))\sigma _2(-1)\sigma _2(1)+4\exp  ({\textstyle{4 \over 3}}\pi\ui\tau (-1))\sigma _1(-1)\sigma _1(1)+4\sigma _2(-1)\sigma _1(1)\cr
  &+16{\exp }({\textstyle{2 \over 3}}\pi\ui \tau (-1))\sigma _2(-1)\sigma _2(1)\sigma _1(-1)\sigma _1(1).\cr}\eqno(5.26d)$$
It turns out that
$$T_{11}T_{22}-T_{12}T_{21}=1.\eqno(5.27)$$
We are looking for a multiplicative solution such that (5.2) holds.
It then follows from (5.1), 5.2), (5.24) that
$$\eqalign{&(T_{11}-p)\alpha +T_{21}\beta =0,\cr
  &T_{12}\alpha +(T_{22}-p)\beta =0\cr}\eqno(5.28)$$
and, as a consequence, that
$$(T_{11}-p)(T_{22}-p)=T_{12}T_{21}\eqno(5.29)$$
or
$$p^2-(T_{11}+T_{22})p+1=0,\eqno(5.30)$$
where (5.27) has been used. The roots $p_1$, $p_2$ of this equation satisfy
$$p_1+p_2=T_{11}+T_{22},\eqno(5.31)$$
and they may be represented in terms of one (not necessarily real) parameter $\omega$  as
$$p_1=\exp (-2\pi\ui\omega ),\quad p_2=\exp (2\pi\ui\omega ).\eqno(5.32)$$
Then we have 
$$p_1+p_2=2\cos (2\pi \omega )=T_{11}+T_{22},\eqno(5.33)$$ 
and the final result is

{\bf Theorem 3.} {\it The characteristic exponent $\omega$ of the multiplicative solutions is given by} 
$$\cos (2\pi \omega )=\cos (2\pi \tau (1))+X,\eqno(5.34)$$
{\it where}
$$\eqalign{&X=2\exp {({\textstyle{4 \over 3}}\pi\ui\tau (-1))}\sigma _0(-1)\sigma _0(1)+2\exp {({\textstyle{4 \over 3}}\pi\ui\tau (1))}\sigma _1(1)\sigma _0(-1)+2\sigma _2(1)\sigma _0(-1)\cr
  &+2\exp {({\textstyle{4 \over 3}}\pi\ui\tau (1))}\sigma _2(1)\sigma _1(-1)+2\exp  {({\textstyle{4 \over 3}}\pi\ui\tau (1))}\sigma _2(-1)\sigma _0(1)+2\sigma _1(-1)\sigma _0(1)\cr
  &+2\exp {({\textstyle{4 \over 3}}\pi\ui\tau (-1))}\sigma _2(-1)\sigma _2(1)+2\exp  {({\textstyle{4 \over 3}}\pi\ui\tau (-1))}\sigma _1(-1)\sigma _1(1)+2\sigma _2(-1)\sigma _1(1)\cr
  &+8\exp ({\textstyle{2 \over 3}}\pi\ui \tau (-1))\sigma _2(-1)\sigma _2(1)\sigma _1(-1)\sigma _1(1)\cr
  &+8\exp ({\textstyle{2 \over 3}}\pi\ui\tau (1))\sigma _2(-1)\sigma _2(1)\sigma _1(-1)\sigma _0(1)+8\exp ({\textstyle{2 \over 3}}\pi\ui\tau (1))\sigma _2(-1)\sigma _1(1)\sigma _0(-1)\sigma _0(1)\cr
  &+8\exp ({\textstyle{2 \over 3}}\pi\ui\tau (1))\sigma _2(1)\sigma _1(-1)\sigma _1(1)\sigma _0(-1)+8\exp ({\textstyle{2 \over 3}}\pi\ui\tau (-1))\sigma _2(-1)\sigma _2(1)\sigma _0(-1)\sigma _0(1)\cr
  &+8\exp ({\textstyle{2 \over 3}}\pi\ui\tau (-1))\sigma _1(-1)\sigma _1(1)\sigma _0(-1)\sigma _0(1)\cr
  &+32\sigma _2(-1)\sigma _2(1)\sigma _1(-1)\sigma _1(1)\sigma _0(-1)\sigma _0(1).\cr}\eqno(5.35)$$
For each of the roots $p_1$ and $p_2$ we may determine the ratio of $\alpha $ and $\beta $ from the upper or lower equation of (5.28). Choosing in each case a convenient normalization, we get the desired multiplicative solutions
$$\eqalign{&f_{p1}(z)=(T_{22}-p_1)f_{\infty 1}(z)-T_{12}f_{\infty 2}(z),\cr
  &f_{p2}(z)=(T_{22}-p_2)f_{\infty 1}(z)-T_{12}f_{\infty 2}(z).\cr}\eqno(5.36)$$
Here the lower equation of (5.28) has been used in both cases so that $T_{11}$, which consists of a considerably longer expression than $T_{22}$, does not appear.

Whenever $D_m=0$ for all the $m=1, 2, \ldots, 6$, the origin is a regular singular point of the differential equation with exponents $\omega =L, -L$. We are not able, however, to see analytically that then the lengthy expression (5.35) reduces to
$X=\cos (2\pi L)-\cos (2\pi \tau (1))$, but this is confirmed in examples of numerical computations, as expected.

\bigskip

\noindent
{\bf 
6.  Asymptotic behaviour of the late coefficients of the formal power series solutions}

According to (4.12), the leading terms of the asymptotic behaviour of the $b$-coefficients for large $m$ are given by
$$\eqalign{&b(\kappa ;m,n_1,n_2)\approx -{1 \over \pi }(2\kappa s_0)^{-m}{{\Gamma (1+\mu (\kappa )+m)} \over {\Gamma (1+m)}}\cr
  &\times \sum\limits_{p_1=0}^{n_1} {\sum\limits_{p_2=0}^{n_2} {e(-\kappa ;p_1,p_2)}}\Gamma (-\mu (-\kappa )-{\textstyle{2 \over 3}}p_1-{\textstyle{1 \over 3}}p_2+m)\cr
  &\times {1 \over {(n_1-p_1)! }}(-2\kappa t_{10})^{n_1-p_1}{1 \over {(n_2-p_2)! }}(-2\kappa t_{20})^{n_2-p_2}.\cr}\eqno(6.1)$$
This result may be used to discuss the asymptotic behaviour of the late $a$-coefficients of the formal solutions. Writing the decomposition (2.17) of the $a$-coefficients in terms of the $b$-coefficients separately for each of three consecutive indices $n=3N, 3N+1, 3N+2$, we have
$$a_{3N}(\kappa )=\sum\limits_{l=0}^{2N} {\sum\limits_{(n_1,n_2)\in J_{3l}} {b(\kappa ;N+l,n_1,n_2)}}=b(\kappa ;N,0,0)+\ldots ,\eqno(6.2)$$
$$a_{3N+1}(\kappa )=\sum\limits_{l=1}^{2N+1} {\sum\limits_{(n_1,n_2)\in J_{3l-1}} {b(\kappa ;N+l,n_1,n_2)}}=b(\kappa ;N+1,1,0)+b(\kappa ;N+1,0,2)+\ldots ,\eqno(6.3)$$
$$a_{3N+2}(\kappa )=\sum\limits_{l=1}^{2N+2} {\sum\limits_{(n_1,n_2)\in J_{3l-2}} {b(\kappa ;N+l,n_1,n_2)}}=b(\kappa ;N+1,0,1)+\ldots ,\eqno(6.4)$$
where
$$J_l=\{(n_1,n_2):2n_1+n_2=l\}\subset \Bbb N_0\times\Bbb N_0,\quad l=0,1,\ldots.  \eqno(6.5)$$
Inserting (6.1) and omitting terms of relative order $N^{-1}$, we get from (6.2)
$$\eqalign{a_{3N}(\kappa )&\approx C_{3N}(\kappa )\sum\limits_{l=0}^{2N} {\sum\limits_{(n_1,n_2)\in J_{3l}} {}\sum\limits_{p_1=0}^{n_1} {\sum\limits_{p_2=0}^{n_2} {e(-\kappa ;p_1,p_2)}}}{{\Gamma (-\mu (-\kappa )-{\textstyle{2 \over 3}}p_1-{\textstyle{1 \over 3}}p_2+l+N)} \over {\Gamma (-\mu (-\kappa )+N)}}\cr
  &\times (2\kappa s_0)^{-l}{1 \over {(n_1-p_1)! }}(-2\kappa t_{10})^{n_1-p_1}{1 \over {(n_2-p_2)! }}(-2\kappa t_{20})^{n_2-p_2},\cr}\eqno(6.6)$$
where
$$\eqalign{C_{3N}(\kappa )&=-{1 \over \pi }{{\Gamma (1+\mu (\kappa )+N)\Gamma (-\mu (-\kappa )+N)} \over {\Gamma (1+N)}}(2\kappa s_0)^{-N}\cr
  &\approx -{1 \over \pi }\Gamma (\mu (\kappa )-\mu (-\kappa )+N)(2\kappa s_0)^{-N}.\cr}\eqno(6.7)$$
Let us now discuss (6.6) in detail: Here only such values of $n_1$ and $n_2$ occur for which
$${\textstyle{2 \over 3}}n_1+{\textstyle{1 \over 3}}n_2=l.$$
If again terms of relative order $N^{-1}$ are omitted, the ratio of the gamma functions is equal to 
$$N^{-(2/3)p_1-(1/3)p_2+l}=N^{(2/3)(n_1-p_1)+(1/3)(n_2-p_2)}.$$
The asymptotic behaviour of (6.6) then becomes
$$\eqalign{&a_{3N}(\kappa )\approx C_{3N}(\kappa )\sum\limits_{l=0}^{2N} {\sum\limits_{(n_1,n_2)\in J_{3l}} {}\sum\limits_{p_1=0}^{n_1} {\sum\limits_{p_2=0}^{n_2} {e(-\kappa ;p_1,p_2)}}}\cr
  &\times {1 \over {(n_1-p_1)! }}[N^{2/3}(2\kappa s_0)^{-2/3}(-2\kappa t_{10})]^{n_1-p_1}{1 \over {(n_2-p_2)! }}[N^{1/3}(2\kappa s_0)^{-1/3}(-2\kappa t_{20})]^{n_2-p_2}.\cr}\eqno(6.8)$$
Here the integer power of $2\kappa s_0$ has been splitted in two factors with fractional powers, the meaning of which is given in terms of the phase factor $\eta$  of (5.10) by 
$$(2\kappa s_0)^{-l}=\left\{ \matrix{(2s_0)^{-(2/3)n_1}(2s_0)^{-(1/3)n_2}\quad{\rm{if}}\quad\kappa =1\hfill\cr
  (2s_0)^{-(2/3)n_1}\eta ^{2n_1}(2s_0)^{-(1/3)n_2}\eta ^{n_2}\quad{\rm{if}}\quad\kappa =-1\hfill\cr} \right..\eqno(6.9)$$
With $n/3$ in place of $N$, the terms in (6.8) look like the terms of the expansion of the exponential function
$$\eqalign{EX_0(n):&=\exp (({\textstyle{1 \over 3}}n)^{2/3}(2\kappa s_0)^{-2/3}(-2\kappa t_{10})+({\textstyle{1 \over 3}}n)^{1/3}(2\kappa s_0)^{-1/3}(-2\kappa t_{20}))\cr
  &=1+({\textstyle{1 \over 3}}n)^{2/3}(2\kappa s_0)^{-2/3}(-2\kappa t_{10})+({\textstyle{1 \over 3}}n)^{1/3}(2\kappa s_0)^{-1/3}(-2\kappa t_{20})+\ldots ,\cr}\eqno(6.10)$$
but because of the three possible values of a third root there are two other such functions,
$$\eqalign{EX_1(n):&=\exp (\eta ^4({\textstyle{1 \over 3}}n)^{2/3}(2\kappa s_0)^{-2/3}(-2\kappa t_{10})+\eta ^2({\textstyle{1 \over 3}}n)^{1/3}(2\kappa s_0)^{-1/3}(-2\kappa t_{20}))\cr
  &=1+\eta ^4({\textstyle{1 \over 3}}n)^{2/3}(2\kappa s_0)^{-2/3}(-2\kappa t_{10})+\eta ^2({\textstyle{1 \over 3}}n)^{1/3}(2\kappa s_0)^{-1/3}(-2\kappa t_{20})+\ldots ,\cr}\eqno(6.11)$$
$$\eqalign{EX_2(n):&=\exp (\eta ^2({\textstyle{1 \over 3}}n)^{2/3}(2\kappa s_0)^{-2/3}(-2\kappa t_{10})+\eta ^4({\textstyle{1 \over 3}}n)^{1/3}(2\kappa s_0)^{-1/3}(-2\kappa t_{20}))\cr
  &=1+\eta ^2({\textstyle{1 \over 3}}n)^{2/3}(2\kappa s_0)^{-2/3}(-2\kappa t_{10})+\eta ^4({\textstyle{1 \over 3}}n)^{1/3}(2\kappa s_0)^{-1/3}(-2\kappa t_{20})+\ldots .\cr}\eqno(6.12)$$
Therefore (6.8) asymptotically shows exponential behaviour given by a certain linear combination of these three exponential functions. They can be identified by the constant term and the two linear terms of their expansion shown above. It is convenient, and easy because of the properties (5.11) of $\eta$, to introduce three functions which are linear combinations of the exponential functions such that only one of the identifying terms is different from zero,
$$\eqalign{&L_0(n):={\textstyle{1 \over 3}}[EX_0(n)+EX_1(n)+EX_2(n)]=1+\ldots ,\cr
  &L_1(n):={\textstyle{1 \over 3}}[EX_0(n)+\eta ^2EX_1(n)+\eta ^4EX_2(n)]=({\textstyle{1 \over 3}}n)^{2/3}(2\kappa s_0)^{-2/3}(-2\kappa t_{10})+\ldots ,\cr
  &L_2(n):={\textstyle{1 \over 3}}[EX_0(n)+\eta ^4EX_1(n)+\eta ^2EX_2(n)]=({\textstyle{1 \over 3}}n)^{1/3}(2\kappa s_0)^{-1/3}(-2\kappa t_{20})+\ldots .\cr}\eqno(6.13)$$
We can now determine the asymptotic behaviour of (6.8) by looking at the contributions from $(p_1,p_2)=(n_1,n_2)$  or  $(p_1,p_2)=(n_1-1,n_2)$ or $(p_1,p_2)=(n_1,n_2-1)$, respectively. This yields
$$a_{3N}(\kappa )\approx C_{3N}(\kappa )\left\{ {S_0(\kappa )L_0(3N)+(2\kappa s_0)^{-1/3}S_1(\kappa )L_1(3N)+(2\kappa s_0)^{-2/3}S_2(\kappa )L_2(3N)} \right\}.\eqno(6.14)$$
In a similar way we may obtain from (6.3) and (6.4), respectively,
$$\eqalign{&a_{3N+1}(\kappa )\approx C_{3N+1}(\kappa )\cr
  &\times \left\{ {(2\kappa s_0)^{-2/3}S_2(\kappa )L_0(3N+1)+S_0(\kappa )L_1(3N+1)+(2\kappa s_0)^{-1/3}S_1(\kappa )L_2(3N+1)} \right\},\cr}\eqno(6.15)$$
$$\eqalign{&a_{3N+2}(\kappa )\approx C_{3N+2}(\kappa )\cr
  &\times \left\{ {(2\kappa s_0)^{-1/3}S_1(\kappa )L_0(3N+2)+(2\kappa s_0)^{-2/3}S_2(\kappa )L_1(3N+2)+S_0(\kappa )L_2(3N+2)} \right\}.\cr}\eqno(6.16)$$
If we now switch back to a representation in terms of the exponential functions (6.10)-(6.12), Stokes multipliers, according to (5.13) above, appear as their factors,
$$\eqalign{&a_{3N}(\kappa )\approx {\textstyle{1 \over 3}}C_{3N}(\kappa )\cr
  &\times \left\{ {\sigma _0(\kappa )EX_0(3N)+\sigma _1(\kappa )EX_1(3N)+\sigma _2(\kappa )EX_2(3N)} \right\},\cr
  &a_{3N+1}(\kappa )\approx {\textstyle{1 \over 3}}C_{3N+1}(\kappa )\cr
  &\times \left\{ {\sigma _0(\kappa )EX_0(3N+1)+\eta ^2\sigma _1(\kappa )EX_1(3N+1)+\eta ^4\sigma _2(\kappa )EX_2(3N+1)} \right\},\cr
  &a_{3N+2}(\kappa )\approx {\textstyle{1 \over 3}}C_{3N+2}(\kappa )\cr
  &\times \left\{ {\sigma _0(\kappa )EX_0(3N+2)+\eta ^4\sigma _1(\kappa )EX_1(3N+2)+\eta ^2\sigma _2(\kappa )EX_2(3N+2)} \right\}.\cr}\eqno(6.17)$$
These three asymptotic equations can be combined to give

{\bf Theorem 4.} {\it The leading terms of the asymptotic exponential behaviour, as $n \to \infty$, of the coefficients of the formal solutions are given by}
$$\eqalign{&a_n(\kappa )\approx -{\textstyle{1 \over 3}}\pi ^{-1}(2\kappa s_0)^{-n/3}\Gamma ({\textstyle{1 \over 3}}[\tau (\kappa )-\tau (-\kappa )+n])\cr
  &\times [\sigma _0(\kappa )EX_0(n)+\eta ^{2n}\sigma _1(\kappa )EX_1(n)+\eta ^{4n}\sigma _2(\kappa )EX_2(n)].\cr}\eqno(6.18)$$
This is a concrete example, with all the quantities determined explicitly, of the structural results obtained by Immink [6], and it is interesting also in the context of related work by other authors [5, 15].
\bigskip

\noindent
{\bf 
7.  Postponed proofs}
\bigskip
\noindent
{\it
7.1  Proof of Lemma 1}

Let
$$f(z)=z^\lambda u(z)\eqno(7.1)$$
and multiply the differential equation for $u(z)$ by $z^6$ in order to remove all the negative powers of $z$. We then are concerned with the differential equation
$$\eqalign{L_zu(z):=&z^8u''+(2\lambda +1)z^7u'\cr
  &-[\sum\limits_{m=1}^6 {D_mz^{6-m}}+(L^2-\lambda ^2)z^6+\sum\limits_{m=1}^6 {B_mz^{6+m}}]u(z)=0.\cr}\eqno(7.2)$$
We are looking for a solution of this differential equation in the form of an integral representation
$$u(z)=\int_{C_{t_2}} {\int_{C_{t_1}} {\int_{C_s} {K(z;s,t_1,t_2)v(s,t_1,t_2)}\d s\d t_1\d t_2}}\eqno(7.3)$$
with the kernel
$$K=K(z;s,t_1,t_2)=\exp (z^3s+z^2t_1+zt_2).\eqno(7.4)$$
If we perform the differentiations with respect to $z$ under the integrals, the differential equation becomes
$$L_zu=\int\!\!\!\int\!\!\!\int {vL_zK}\d s\d t_1\d t_2=0,\eqno(7.5)$$
where
$$\eqalign{{1 \over K}L_zK&=(9s^2-B_6)z^{12}+(12st_1-B_5)z^{11}+(6st_2+4t_1^2-B_4)z^{10}\cr
  &+([6\lambda +9]s+4t_1t_2-B_3)z^9+([4\lambda +4]t_1+t_2^2-B_2)z^8\cr
  &+([2\lambda +1]t_2-B_1)z^7+(\lambda ^2-L^2)z^6-\sum\limits_{m=1}^6 {D_mz^{6-m}}.\cr}\eqno(7.6)$$
It is advisable to rewrite this in terms of the quantities $p_1$, $p_2$, $p_3$, according to (1.5), which determine the exponential factors of the formal solutions (1.2)-(1.3), or in terms of the related quantities  $s_0$, $t_{10}$, $t_{20}$  according to (2.3).  Then (7.6) becomes
$$\eqalign{{1 \over K}L_zK&=9(s^2-s_0^2)z^{12}+12(st_1-s_0t_{10})z^{11}+(6[st_2-s_0t_{20}]+4[t_1^2-t_{10}^2])z^{10}\cr
  &+([6\lambda +9]s+4[t_1t_2-t_{10}t_{20}]+[4t_{10}t_{20}-B_3])z^9\cr
  &+([4\lambda +4]t_1+[t_2^2-t_{20}^2]+[t_{20}^2-B_2])z^8\cr
  &+([2\lambda +1]t_2-B_1)z^7+(\lambda ^2-L^2)z^6-\sum\limits_{m=1}^6 {D_mz^{6-m}}.\cr}\eqno(7.7)$$
 By repeated partial integrations of the exponential function with respect to $s$ or $t_1$ or $t_2$, respectively, it is possible to get rid of all the powers of $z$. This lengthy task can more conveniently be performed in a formal way as follows: 
Let us find a partial differential expression $M=M_{s,t_1,t_2}$ with respect to $s$, $t_1$, $t_2$, independent of $z$, such that
$$L_zK\equiv M_{s,t_1,t_2}K.\eqno(7.8)$$
The powers of $z$ correspond to partial derivatives with respect to $s$ or $t_1$ or $t_2$, and there seems to be some ambiguity as to the choice of the partial derivatives which yield the same total power of $z$. This ambiguity is resolved by the following reasonning: In order to get the order of each derivative as small as possible, we would prefer derivatives with respect to $s$, which account for $z^3$, with highest priority, next  $t_1$, which accounts for $z^2$, last $t_2$, which accounts only for $z$. Conflicting with this policy, however, are some other requirements which ensure that we get a differential expression appropriate for our purpose. So the singularities should be at the right places, which are related to the coefficients of the exponential factor of the formal solutions we want to represent. In particular, bilinear or quadratic factors, such as $t_1^2$ for instance, should occur in the form $t_1^2-t_{10}^2$. For this reason we have in (7.7) already added and subtracted the term $4t_{10}^2$. It is then necessary to treat these two terms differently, so that any power is always multiplied by the corresponding derivative. In the example just mentioned this means that $z^{10}$ has to be translated into $\partial ^4/\partial s^3\partial t_2$ for one term, but into $\partial ^4/\partial s^2 \partial t_1^2$ for the other. The unique result now is
$$\eqalign{M_{s,t_1,t_2}&=9(s^2-s_0^2)\partial ^4/\partial s^4+12(st_1-s_0t_{10})\partial ^4/\partial s^3\partial t_1\cr
  &+6(st_2-s_0t_{20})\partial ^4/\partial s^3\partial t_2+4(t_1^2-t_{10}^2)\partial ^4/\partial s^2\partial t_1^2\cr
  &+([6\lambda +9]s+[4t_{10}t_{20}-B_3])\partial ^3/\partial s^3+4[t_1t_2-t_{10}t_{20}]\partial ^4/\partial s^2\partial t_1\partial t_2\cr
  &+([4\lambda +4]t_1+[t_{20}^2-B_2])\partial ^3/\partial s^2\partial t_1+[t_2^2-t_{20}^2]\partial ^4/\partial s^2\partial t_2^2\cr
  &+([2\lambda +1]t_2-B_1)\partial ^3/\partial s^2\partial t_2+(\lambda ^2-L^2)\partial ^2/\partial s^2\cr
  &-D_1\partial ^2/\partial s\partial t_1-D_2\partial ^2/\partial s\partial t_2-D_3\partial /\partial s-D_4\partial /\partial t_1-D_5\partial /\partial t_2-D_6.\cr}\eqno(7.9)$$
Next let us introduce the adjoint differential expression $\overline M=\overline M_{s,t_1,t_2}$ defined, with any sufficiently differentiable function $v=v(s,t_1,t_2)$, by
$$\eqalign{&\overline M_{s,t_1,t_2}v=(\partial ^4/\partial s^4)\{9(s^2-s_0^2)v\}+(\partial ^4/\partial s^3\partial t_1)\{12(st_1-s_0t_{10})v\}\cr
  &+(\partial ^4/\partial s^3\partial t_2)\{6(st_2-s_0t_{20})v\}+(\partial ^4/\partial s^2\partial t_1^2)\{4(t_1^2-t_{10}^2)v\}\cr
  &-(\partial ^3/\partial s^3)\{([6\lambda +9]s+[4t_{10}t_{20}-B_3])v\}+(\partial ^4/\partial s^2\partial t_1\partial t_2)\{4(t_1t_2-t_{10}t_{20})v\}\cr
  &-(\partial ^3/\partial s^2\partial t_1)\{([4\lambda +4]t_1+[t_{20}^2-B_2])v\}+(\partial ^4/\partial s^2\partial t_2^2)\{(t_2^2-t_{20}^2)v\}\cr
  &-(\partial ^3/\partial s^2\partial t_2)\{([2\lambda +1]t_2-B_1)v\}+(\lambda ^2-L^2)(\partial ^2v/\partial s^2)\cr
  &-D_1(\partial ^2v/\partial s\partial t_1)-D_2(\partial ^2v/\partial s\partial t_2)+D_3(\partial v/\partial s)+D_4(\partial v/\partial t_1)+D_5(\partial v/\partial t_2)-D_6v.\cr}\eqno(7.10)$$
The difference as compared with $M$ is that the factors in front of each derivative  have here to be differentiated too and that all the terms of odd order change their sign.  
The usefulness of the adjoint expression lies in the formula, known in the one-variable case as the identity of Lagrange, 
$$vMK-K\overline Mv=RHS,\eqno(7.11)$$
where the right-hand side $RHS$ is a lengthy expression, bilinear in $v$ and $K$ or their partial derivatives, which may be so arranged that it consists of a sum of terms of which each is a total derivative with respect to one of the variables. For the term with the factor $-D_1$, for instance, this reads
$$v{{\partial ^2K} \over {\partial s\partial t_1}}-K{{\partial ^2v} \over {\partial s\partial t_1}}={\textstyle{1 \over 2}}{\partial  \over {\partial s}}[v{{\partial K} \over {\partial t_1}}-K{{\partial v} \over {\partial t_1}}]+{\textstyle{1 \over 2}}{\partial  \over {\partial t_1}}[v{{\partial K} \over {\partial s}}-K{{\partial v} \over {\partial s}}].\eqno(7.12)$$
By means of (7.5) and (7.8), the equation to be satisfied now is
$$\eqalign{L_zu=&\int\!\!\!\int\!\!\!\int {vM_{s,t_1,t_2}K\d s\d t_1\d t_2}\cr
  =&\int\!\!\!\int\!\!\!\int {K\overline M_{s,t_1,t_2}v\d s\d t_1\d
t_2}+\int\!\!\!\int\!\!\!\int {RHS\d s\d t_1\d t_2}=0.\cr}\eqno(7.13)$$
The first term in the second line can be made to vanish if $v$ is required to be a solution of the partial differential equation
$$\overline M_{s,t_1,t_2}v(s,t_1,t_2)=0,\eqno(7.14)$$
and the second term by a suitable choice of the contours of integration, for it is a sum of semi-integrated terms, each involving the difference of the values of the integrand at the termini of the contour of one variable and only two remaining integrals with respect to the other two variables. The partial differential equation (2.2) is the same as (7.14), after the derivatives of the products in (7.10) have been resolved. This completes the proof of Lemma 1.
\bigskip

\noindent
{\it 
7.2  Proof of Lemma 2}

In terms of the shifted variables
$$S=s-\kappa s_0,\quad T_1=t_1-\kappa t_{10},\quad T_2=t_2-\kappa t_{20},\eqno(7.15)$$
the differential equation (2.2) reads
$$\eqalign{&9S(S+2\kappa s_0)(\partial ^4v/\partial S^4)+12(ST_1+\kappa t_{10}S+\kappa s_0T_1)(\partial ^4v/\partial S^3\partial T_1)\cr
  &+6(ST_2+\kappa t_{20}S+\kappa s_0T_2)(\partial ^4v/\partial S^3\partial T_2)+4T_1(T_1+2\kappa t_{10})(\partial ^4v/\partial S^2\partial T_1^2)\cr
  &+T_2(T_2+2\kappa t_{20})(\partial ^4v/\partial S^2\partial T_2^2)+4(T_1T_2+\kappa t_{20}T_1+\kappa t_{10}T_2)(\partial ^4v/\partial S^2\partial T_1\partial T_2)\cr
  &+([81-6\lambda ]\kappa s_0+B_3-4t_{10}t_{20})(\partial ^3v/\partial S^3)+(81-6\lambda )S(\partial ^3v/\partial S^3)\cr
  &+([52-4\lambda ]\kappa t_{10}+B_2-t_{20}^2)(\partial ^3v/\partial S^2\partial T_1)+(52-4\lambda )T_1(\partial ^3v/\partial S^2\partial T_1)\cr
  &+([25-2\lambda ]\kappa t_{20}+B_1)(\partial ^3v/\partial S^2\partial T_2)+(25-2\lambda )T_2(\partial ^3v/\partial S^2\partial T_2)\cr
  &+([\lambda -12]^2-L^2)(\partial ^2v/\partial S^2)-D_1(\partial ^2v/\partial S\partial T_1)-D_2(\partial ^2v/\partial S\partial T_2)\cr
  &+D_3(\partial v/\partial S)+D_4(\partial v/\partial T_1)+D_5(\partial v/\partial T_2)-D_6v=0.\cr}\eqno(7.16)$$
Inserting a power series solution
$$v=\sum\limits_{m=0} {\sum\limits_{n_1=0} {\sum\limits_{n_2=0} {a(m,n_1,n_2)S^{\mu +m}T_1^{-\nu _1-n_1}T_2^{-\nu _2-n_2}}}},\eqno(7.17)$$
we obtain for the coefficients the recurrence relation
$$\eqalign{&6\kappa s_0(\mu +m)(\mu +m-1)(\mu +m-2)[3(\mu +m)-2(\nu _1+n_1)-(\nu _2+n_2)-\lambda -\tau (\kappa )+6]\cr
  &\times a(m,n_1,n_2)\cr
  &+(\mu +m-1)(\mu +m-2)\{[3(\mu +m)-2(\nu _1+n_1)-(\nu _2+n_2)-\lambda +3]^2-L^2\}\cr
  &\times a(m-1,n_1,n_2)\cr
  &+(\mu +m-1)(\mu +m-2)(\nu _1+n_1-1)\{4\kappa t_{10}[-3(\mu +m)+2(\nu _1+n_1)+(\nu _2+n_2)\cr  
  &+\lambda -4]
  +t_{20}^2-B_2\}a(m-1,n_1-1,n_2)\cr
  &+(\mu +m-1)(\mu +m-2)(\nu _2+n_2-1)\{2\kappa t_{20}[-3(\mu +m)+2(\nu _1+n_1)+(\nu _2+n_2)\cr
  &+\lambda -{\textstyle{7 \over 2}}]
  +B_1\}a(m-1,n_1,n_2-1)\cr
  &+D_1(\mu +m-2)(\nu _1+n_1-1)a(m-2,n_1-1,n_2)\cr
  &+D_2(\mu +m-2)(\nu _2+n_2-1)a(m-2,n_1,n_2-1)\cr
  &+D_3(\mu +m-2)a(m-2,n_1,n_2)-D_4(\nu _1+n_1-1)a(m-3,n_1-1,n_2)\cr
  &-D_5(\nu _2+n_2-1)a(m-3,n_1,n_2-1)-D_6a(m-3,n_1,n_2)=0,\cr}\eqno(7.18)$$
valid for $m\ge 0$, $n_1\ge 0$, $n_2\ge 0$ provided that we agree that all the $a$-coefficients are equal to zero if any of the indices $m, n_1, n_2$ is less than zero.
Assuming that $a(0,0,0)\ne 0$, we get from the equation for $m=n_1=n_2=0$ the indicial equation
$$\mu (\mu -1)(\mu -2)(3\mu -2\nu _1-\nu _2-\lambda -\tau (\kappa )+6)=0,\eqno(7.19)$$
with $\tau (\kappa )$ according to (1.6). Possible values of the exponent $\mu $ are therefore $\mu =0, 1, 2, $ or
$$\mu ={\textstyle{1 \over 3}}[2\nu _1+\nu _2+\lambda +\tau (\kappa )]-2.\eqno(7.20)$$
In order to avoid complications, we may assume that the last possibility does not yield an integer value. This can always be guarantied be a suitable choice of the still disposable parameter $\lambda $. Inspection of the recurrence relation then shows that each of the possible values of $\mu $ leads to a solution of the partial differential equation. We here need not further consider the three solutions which are regular with respect to $S$ at $S=0$, but we continue to discuss the singular one with the exponent (7.20), writing for the associated coefficients $A(\kappa ; m, n_1, n_2)$ rather than $a(m, n_1, n_2)$. Simplifying by means of (7.20) and introducing the $b$-coefficients according to (2.6) , we get (2.7). The coefficients for which $3m-2n_1-n_2=0$ are constants of integration and may be chosen to be zero. This completes the proof of Lemma 2.
\bigskip

\noindent
{\it
7.3  Proof of the limit formula}

We have to verify that the limit formula (4.14) satisfies (4.13 ). Substituting it for the $e$-coefficients in (4.13) and interchanging the summations, we have to evaluate sums such as
$$\sum\limits_{p=j}^n {{{(-1)^{n-p}} \over {(p-j)! (n-p) !}}}={{(-1)^{n-j}} \over {(n-j)! }}\sum\limits_{q=0}^{n-j} {(-1)^q{{n-j} \choose q}
 }=\left\{ \matrix{1\quad  {\rm{if}}\quad j=n\hfill\cr
  0\quad  {\rm{if}}\quad 0\le j\le n-1,\;n>0\hfill\cr} \right..\eqno(7.21)$$
Therefore only one term survives on the right-hand side, which then becomes equal to the left.
\bigskip

\noindent
{\it 
7.4  Proof of Theorem 1}

In order to prove Theorem 1, we first multiply the continuation formula (3.5), with $r_1=r_2=0$, by
$$({\textstyle{1 \over 2}}-{s \over {2\kappa s_0}})^{-\mu (\kappa )}=\sum\limits_{j=0}^\infty  {{{(\mu (\kappa ))_j} \over {j! }}}({\textstyle{1 \over 2}}+{s \over {2\kappa s_0}})^j,\eqno(7.22)$$
where the left-hand side is used on the left and the right-hand side on the right of the continuation formula. Then the left-hand side is the same as above (4.1) with the power factor in front of the series removed, but the right becomes, after multiplication of the two power series, 
$$\eqalign{&\sum\limits_{q_1=0}^\infty  {\sum\limits_{q_2=0}^\infty  {E(-\kappa ;}}0,0;q_1,q_2)\Gamma (-\mu (-\kappa )-{\textstyle{2 \over 3}}q_1-{\textstyle{1 \over 3}}q_2)\cr
  &\times \sum\limits_{k=0}^\infty  {\sum\limits_{l_1=0}^k {\sum\limits_{l_2=0}^k {\Gamma (\nu _1+q_1+l_1)\Gamma (\nu _2+q_2+l_2)H(-\kappa ;k,l_1,l_2;q_1,q_2)}}}\cr
  &\times ({\textstyle{1 \over 2}}+{s \over {2\kappa s_0}})^{\mu (-\kappa )+{\textstyle{2 \over 3}}q_1+{\textstyle{1 \over 3}}q_2+k}(t_1+\kappa t_{10})^{-\nu _1-q_1-l_1}(t_2+\kappa t_{20})^{-\nu _2-q_2-l_2},\cr}\eqno(7.23)$$
with $H$ defined in (4.16).  Proceeding as above, and making use of the formula
$$(\alpha -k)_m=(\alpha )_m{{(1-\alpha )_k} \over {(1-\alpha -m)_k}},\eqno(7.24)$$
we obtain
$$\eqalign{&{{-\pi } \over {\sin (\pi \mu (\kappa ))}}{1 \over {\Gamma (1+\mu (\kappa )+m)}}(2\kappa s_0)^mb(\kappa ;m,n_1,n_2)\cr
  &\sim \sum\limits_{q_1=0}^{n_1} {\sum\limits_{q_2=0}^{n_2} {E(-\kappa ;}}0,0;q_1,q_2){{\Gamma (-\mu (-\kappa )-{\textstyle{2 \over 3}}q_1-{\textstyle{1 \over 3}}q_2+m)} \over {m! }}\cr
  &\times \sum\limits_{k=0}^\infty  {\sum\limits_{l_1=0}^k {\sum\limits_{l_2=0}^k {{{(1+\mu (-\kappa )+{\textstyle{2 \over 3}}q_1+{\textstyle{1 \over 3}}q_2)_k} \over {(1+\mu (-\kappa )+{\textstyle{2 \over 3}}q_1+{\textstyle{1 \over 3}}q_2-m)_k}}H(-\kappa ;k,l_1,l_2;q_1,q_2)}}}\cr
  &\times \left( {{1 \over {j_1! }}(-2\kappa t_{10})^{j_1}} \right)_{j_1+q_1+l_1=n_1}\left( {{{(1} \over {j_2! }}(-2\kappa t_{20})^{j_2}} \right)_{j_2+q_2+l_2=n_2}.\cr}\eqno(7.25)$$
 Keeping the first $K+1$ singular terms on the right and solving for the $E$-coefficient with $q_1=n_1$, $q_2=n_2$, we have
$$\eqalign{&E(-\kappa ;0,0;n_1,n_2)\cr
  &\times [1+\sum\limits_{k=1}^K {{{(1+\mu (-\kappa )+{\textstyle{2 \over 3}}n_1+{\textstyle{1 \over 3}}n_2)_k} \over {(1+\mu (-\kappa )+{\textstyle{2 \over 3}}n_1+{\textstyle{1 \over 3}}n_2-m)_k}}H(-\kappa ;k,0,0;n_1,n_2)}+O(m^{-K-1})]\cr
  &={{-\pi } \over {\sin (\pi \mu (\kappa ))}}{{m! } \over {\Gamma (1+\mu (\kappa )+m)\Gamma (-\mu (-\kappa )-{\textstyle{2 \over 3}}n_1-{\textstyle{1 \over 3}}n_2+m)}}(2\kappa s_0)^mb(\kappa ;m,n_1,n_2)\cr
  &-\mathop {\sum\limits_{q_1=0}^{n_1} {\sum\limits_{q_2=0}^{n_2} {}}}\limits_{(q_1,q_2)\ne (n_1,n_2)}E(-\kappa ;0,0;q_1,q_2){{\Gamma (-\mu (-\kappa )-{\textstyle{2 \over 3}}q_1-{\textstyle{1 \over 3}}q_2+m)} \over {\Gamma (-\mu (-\kappa )-{\textstyle{2 \over 3}}n_1-{\textstyle{1 \over 3}}n_2+m)}}\cr
  &\times [1+\sum\limits_{k=1}^K {\sum\limits_{l_1=0}^k {\sum\limits_{l_2=0}^k {{{(1+\mu (-\kappa )+{\textstyle{2 \over 3}}q_1+{\textstyle{1 \over 3}}q_2)_k} \over {(1+\mu (-\kappa )+{\textstyle{2 \over 3}}q_1+{\textstyle{1 \over 3}}q_2-m)_k}}H(-\kappa ;k,l_1,l_2;q_1,q_2)}}}+O(m^{-K-1})]\cr
  &\times \left( {{1 \over {j_1! }}(-2\kappa t_{10})^{j_1}} \right)_{j_1+q_1+l_1=n_1}\left( {{1 \over {j_2! }}(-2\kappa t_{20})^{j_2}} \right)_{j_2+q_2+l_2=n_2}.\cr}\eqno(7.26)$$
This is essentially (4.15) because of (4.11) and completes the proof of Theorem 1.
\bigskip

\noindent
{\it 
7.5  Choice of the computational parameter $\lambda$}

If $D_3=D_6=0$, the recurrence relation (2.7) for $b(\kappa; m, 0, 0)$ reduces to a two-term relation, and we obtain
$$b(\kappa ;m,0,0)=(2\kappa s_0)^{-m}{{(-{\textstyle{1 \over 3}}L+{\textstyle{1 \over 3}}\tau (\kappa ))_m({\textstyle{1 \over 3}}L+{\textstyle{1 \over 3}}\tau (\kappa ))_m} \over {m! }}.\eqno(7.27)$$
Rewriting (4.16) by means of the identity
$${{(x)_{k-j}} \over {(k-j)! }}={{(x)_k} \over {k! }}{{(-k)_j} \over {(1-x-k)_j}},\eqno(7.28)$$
we then have
$$H(-\kappa ;k,0,0;0,0)={{(\mu (\kappa ))_k} \over {k! }}\sum\limits_{j=0}^k {{{(-k)_j(-{\textstyle{1 \over 3}}L+{\textstyle{1 \over 3}}\tau (-\kappa ))_j({\textstyle{1 \over 3}}L+{\textstyle{1 \over 3}}\tau (-\kappa ))_j} \over {(1-\mu (\kappa )-k)_j(1+\mu (-\kappa ))_jj! }}}.\eqno(7.29)$$
Because of 
$3\mu (\kappa )=\lambda +\tau (\kappa )-3$ according to (2.5) and (2.11) and
$\tau (\kappa )=3-\tau (-\kappa )$ according to (1.6), the series is a terminating one-balanced hypergeometric series at unit argument, which can be summed by the theorem of Saalsch\"utz [11], so that
$$H(-\kappa ;k,0,0;0,0)={{(-{\textstyle{1 \over 3}}L+{\textstyle{1 \over 3}}\lambda )_k({\textstyle{1 \over 3}}L+{\textstyle{1 \over 3}}\lambda )_k} \over {(1+\mu (-\kappa ))_kk! }}.\eqno(7.30)$$
In total, we get
$$\eqalign{&e(-\kappa ;0,0)={{-\pi } \over {\Gamma (-{\textstyle{1 \over 3}}L+{\textstyle{1 \over 3}}\tau (\kappa ))\Gamma ({\textstyle{1 \over 3}}L+{\textstyle{1 \over 3}}\tau (\kappa ))}}{{\Gamma (-{\textstyle{1 \over 3}}L+{\textstyle{1 \over 3}}\tau (\kappa )+m)\Gamma ({\textstyle{1 \over 3}}L+{\textstyle{1 \over 3}}\tau (\kappa )+m)} \over {\Gamma ({\textstyle{1 \over 3}}\lambda +{\textstyle{1 \over 3}}\tau (\kappa )+m)\Gamma (-{\textstyle{1 \over 3}}\lambda +{\textstyle{1 \over 3}}\tau (\kappa )+m)}}\cr
  &\times \left[ {1+\sum\limits_{k=1}^K {{{(-{\textstyle{1 \over 3}}L+{\textstyle{1 \over 3}}\lambda )_k({\textstyle{1 \over 3}}L+{\textstyle{1 \over 3}}\lambda )_k} \over {({\textstyle{1 \over 3}}\lambda +{\textstyle{1 \over 3}}\tau (-\kappa )-m)_kk! }}}+O(m^{-K-1})} \right]^{-1}.\cr}\eqno(7.31)$$
For $\lambda=-L$ or  $\lambda=L$, the terms with $k=1, 2, \ldots$ of the asymptotic series all vanish and the factor in front of the series becomes independent of $m$. More generally this means that some terms which might become quite large when L is not small can be removed by such a choice of $\lambda$ from the asymptotic series and incorporated in the $m$-dependence of the function in front of the series.
\vfill\eject
\noindent
{\bf References}

[1]  W. Balser,  W. B. Jurkat, and D. A. Lutz, Transfer of connection problems for meromorphic differential equations of rank $\ge 2$ and representations of solutions, {\it J. Math. Anal. Appl.} {\bf 85} (1982) 488-542.

[2]	B. L. J. Braaksma, Multisummability and Stokes multipliers of linear meromorphic 	differential equations, {\it J. Differ. Equations} {\bf 92} (1991) 45-75.

[3]	W. B\"uhring, The characteristic exponent of second-order linear differential equations with two irregular singular points, {\it Proc. Amer. Math. Soc.} {\bf 118} (1993) 801-812.

[4]  A. Duval, Triconfluent Heun's equation, in: A. Ronveaux, ed., {\it Heun's differential equations} (Oxford University Press, Oxford, New York, Tokyo, 1995).

[5]  R. Hoeppner and R. Sch\"afke, On the remainders of asymptotic expansions of solutions of linear differential equations near irregular singular points of higher rank, {\it Math. Nachr.} {\bf 205} (1999) 89-113.

[6]	G. K. Immink, A note on the relationship between Stokes multipliers and formal solutions of 	analytic differential equations, {\it SIAM J. Math. Anal.} {\bf 21} (1990) 782-792.

[7]	E. L. Ince, {\it Ordinary differential equations} ( Dover, New York, 1956).

[8]	M. Kohno, Stokes phenomenon for general linear ordinary differential equations with two 	singular points, in: Japan - United States seminar on ordinary differential and functional equations, {\it Lecture Notes in Mathematics} {\bf 243} (Springer, 1971)  310-314. 

[9]	M. Kohno,  A two point connection problem for general linear ordinary differential equations, {\it Hiroshima Math. J.} {\bf 4} (1974) 293-338.
 	
[10]  M. Loday-Richaud, Calcul des invariants de Birkhoff des syst\`emes d'ordre deux, {\it Funk. Ekvac.} {\bf 33} (1990) 161-225.

[11]  Y. L. Luke,  {\it The special functions and their approximations, Vol. 1} (Academic Press,  New York, 1969).

[12]	J. Martinet and J.-P. Ramis, Elementary acceleration and multisummability I, {\it Ann. Inst. Henri Poincar\'e, Phys. Th\'eor.} {\bf 54} (1991) 331-401.

[13]	F. Naundorf,  A connection problem for second order linear differential equations with two irregular singular points, {\it SIAM J. Math. Anal.} {\bf 7} (1976) 157-175.

[14]	F. Naundorf,  Ein Verfahren zur Berechnung der charakteristischen Exponenten von linearen 	Differentialgleichungen zweiter Ordnung mit zwei stark singul\"aren Stellen, {\it Z. Angew. Math. Mech.} {\bf 57} (1977) 47-49.

[15]  A. B. Olde Daalhuis and F. W. J. Olver, On the calculation of Stokes multipliers for linear differential equations of the second order, {\it Methods Appl. Anal.} {\bf 2} (1995) 348-367.

[16]	F. W. J. Olver,  {\it Asymptotics and special functions}  (Academic Press, New York 1974).

[17]	R. Sch\"afke and D. Schmidt,  The connection problem for general linear ordinary differential 	equations at two regular singular points with applications in the theory of special functions,  {\it	SIAM J. Math. Anal.} {\bf 11} (1980) 848-862.

[18]  Y. Sibuya, {\it Global theory of a second-order linear ordinary differential equation with a polynomial coefficient}  (North-Holland, Amsterdam, New York, 1975).

[19]  V. S. Varadarajan, Linear meromorphic differential equations: a modern point of view, {\it  Bull. Amer. Math. Soc.} {\bf 33} (1996) 1-42.

[20]	E. Wagenf\"uhrer,  Die  Determinantenmethode zur Berechnung des charakteristischen 	Exponenten der endlichen Hillschen Differentialgleichung, {\it Numer. Math.} {\bf 35} (1980) 405-420.

[21]	E. Wagenf\"uhrer and H. Lang,  Berechnung des charakteristischen Exponenten der endlichen Hillschen Differentialgleichung durch numerische Integration,  {\it Numer. Math.} {\bf 32} (1979) 	31-50.

[22]  E. T. Whittaker and G. N. Watson, {\it A course of modern analysis} (Cambridge University Press, London , 1927).
\bigskip

Wolfgang B\"uhring

Physikalisches Institut

Universit\"at Heidelberg

Philosophenweg 12

69120 Heidelberg

GERMANY
\bigskip

buehring@physi.uni-heidelberg.de

\bigskip

\vfill\eject
\bye